\documentclass[]{amsart}
\numberwithin{equation}{section}
\usepackage{amssymb}
\newcommand{\lbl}[1]{{\tt \small [#1]}\label{#1}}
\renewcommand{\lbl}[1]{\label{#1}}

\newcommand{\xxx}{\mbox{\fontfamily{phv}\selectfont x}}
\newcommand{\yyy}{\mbox{\fontfamily{phv}\selectfont y}}
\newcommand{\D}{\mbox{\fontfamily{phv}\selectfont D}_q}
\newcommand{\Z}{\mbox{\fontfamily{phv}\selectfont Z}}
\newcommand{\X}{\mbox{\fontfamily{phv}\selectfont X}}

\newcommand{\I}{\mbox{\fontfamily{phv}\selectfont I}}

\def\<{\langle}
\def\>{\rangle}

\def\<{\langle}
\def\>{\rangle}

\def\v{{\bf v}}

\def\w{{\bf w}}

\newcommand{\rf}[1]{(\ref{#1})}

\newcommand{\cal}{\mathcal}

\newcommand{\calF}{{\mathcal F}}

\newcommand{\aaa}{{\mathbf a}}
\newcommand{\bbb}{{\mathbf b}}

\newcommand{\AAA}{{\mathbf A}}
\newcommand{\BBB}{{\mathbf B}}
\newcommand{\CCC}{{\mathbf C}}
\newcommand{\DDD}{{\mathbf D}}

\newcommand{\II}{{\mathbf I}} 

\newcommand{\sR}{{\mathbb R}}

\newcommand{\sN}{{\mathbb N}}

\newcommand{\mC}{\mathbf{C}}

\newcommand{\eps}{\varepsilon}
\newcommand{\la}{\lambda}

\newcommand{\be}{\begin{equation}}
\newcommand{\ee}{\end{equation}}

      \newtheorem{theorem}{Theorem}[section]
       \newtheorem{proposition}[theorem]{Proposition}
       \newtheorem{corollary}[theorem]{Corollary}
       \newtheorem{lemma}[theorem]{Lemma}
\newtheorem{claim}[theorem]{Claim}
\theoremstyle{remark}
       \newtheorem{remark}{Remark}[section]
\theoremstyle{definition}
\newtheorem{definition}{Definition}[section]

\newtheorem{example}[theorem]{Example}

\def\<{\langle}
\def\>{\rangle}

\def\<{\langle}
\def\>{\rangle}

\def\v{{\bf v}}

\def\w{{\bf w}}

\renewcommand{\AAA}{{A}_{t,s,u}}
\renewcommand{\BBB}{{B}_{t,s,u}}
\renewcommand{\CCC}{{C}_{t,s,u}}
\renewcommand{\DDD}{{D}_{t,s,u}}
\newcommand{\EEE}{{E}_{t,s,u}}
\newcommand{\FFF}{{F}_{t,s,u}}
\renewcommand{\aaa}{{a}_{t,s,u}}
\renewcommand{\bbb}{{b}_{t,s,u}}

\newcommand{\wo}{\mathbb{E}}                                                    
\newcommand{\wwo}[2]{\mathbb{E}\left[\left.{#1}\right|{#2}\right]}              
\DeclareMathOperator{\Var}{Var}                                                 
\newcommand{\wVar}[2]{\Var\left[\left.{#1}\right|{#2}\right]}                   

\newcommand{\pln}[2]{p_{#1}\left(X_{#2};#2\right)}

\newcommand{\Fsu}{{\mathcal{F}}_{s,u}}

\newcommand{\Fru}{{\mathcal{F}}_{r,u}}
\newcommand{\Frt}{{\mathcal{F}}_{r,t}}

\newcommand{\asru}{a_{s,r,u}}
\newcommand{\bsru}{b_{s,r,u}}
\newcommand{\Asru}{A_{s,r,u}}
\newcommand{\Bsru}{B_{s,r,u}}
\newcommand{\Csru}{C_{s,r,u}}
\newcommand{\Dsru}{D_{s,r,u}}
\newcommand{\Esru}{E_{s,r,u}}
\newcommand{\Fsru}{F_{s,r,u}}

\newcommand{\atru}{a_{t,r,u}}
\newcommand{\btru}{b_{t,r,u}}
\newcommand{\Atru}{A_{t,r,u}}
\newcommand{\Btru}{B_{t,r,u}}
\newcommand{\Ctru}{C_{t,r,u}}
\newcommand{\Dtru}{D_{t,r,u}}
\newcommand{\Etru}{E_{t,r,u}}
\newcommand{\Ftru}{F_{t,r,u}}

\newcommand{\atsu}{a_{t,s,u}}
\newcommand{\btsu}{b_{t,s,u}}
\newcommand{\Atsu}{A_{t,s,u}}
\newcommand{\Btsu}{B_{t,s,u}}
\newcommand{\Ctsu}{C_{t,s,u}}
\newcommand{\Dtsu}{D_{t,s,u}}
\newcommand{\Etsu}{E_{t,s,u}}
\newcommand{\Ftsu}{F_{t,s,u}}

\newcommand{\asrt}{a_{s,r,t}}
\newcommand{\bsrt}{b_{s,r,t}}

\author{
W{\l}odzimierz  Bryc
}
\thanks{\noindent Research partially supported by NSF
grant \#INT-0332062 and by the C.P. Taft Memorial Fund.}
\address{
Department of Mathematics,
University of Cincinnati,
PO Box 210025,
Cincinnati, OH 45221--0025, USA}
\email{Wlodzimierz.Bryc@UC.edu}

\author{Wojciech Matysiak}
\address{Faculty of Mathematics and Information Science\\
Warsaw University of Technology\\
pl. Politechniki 1\\
 00-661 Warszawa, Poland}
\email{matysiak@mini.pw.edu.pl}

\author{Jacek Weso{\l}owski}
\address{ Faculty of Mathematics and Information Science\\
Warsaw University of Technology\\ pl. Politechniki 1\\ 00-661
Warszawa, Poland}
\email{wesolo@alpha.mini.pw.edu.pl}

\date{January 7, 2005; Revised: April 3, 2005}


\title[Quadratic Harnesses]{Quadratic Harnesses,
 $q$-commutations,
and orthogonal martingale polynomials}
\begin{document}
\begin{abstract} We introduce the quadratic harness condition and show that integrable quadratic harnesses have orthogonal
martingale polynomials with a three step recurrence that satisfies a $q$-commutation relation.
 This
implies that  quadratic harnesses are essentially determined
 uniquely by five numerical constants. Explicit recurrences
for the orthogonal martingale polynomials are derived in several cases of interest.
\end{abstract}


\maketitle
%

\tolerance=2000

\section{Introduction}
Hammersley \cite{Hammersley} introduced harnesses on $\sR^n$ as the probabilistic models of  long-range
misorientation in the crystalline structure of metals. Several authors studied mathematical aspects of
the concept: Mansuy and Yor \cite{Mansuy-Yor04} analyze  harnesses on $\sR_+$, Williams
\cite{Williams} analyzes  harnesses with arbitrary discrete index set; see also
\cite{MR82h:60152}, \cite{MR999888}, \cite{MR93h:60129}, \cite{MR91e:60158}.

The  class of random fields on $\sR_+$
 which we call the
quadratic harnesses is related to Hammersley's harnesses
by mimicking the relation between the martingale and the quadratic martingale conditions.
Such processes have already been studied implicitly by several authors, see the paragraph following Definition \ref{Def QH}.
In particular, ref. \cite{Bryc-Wesolowski-03} gives a construction of the three parameter family of Markov processes
with the quadratic harness property.  Some of these processes are related to the free L\'evy processes,
and some correspond to the non-commutative $q$-Gaussian processes introduced by
Frisch and Bourret \cite{FB70} and studied in \cite{BKS97}.

In full generality,  the constructions of the quadratic harnesses on $\sR_+$
are not yet well understood. In this paper we  concentrate on
uniqueness without dealing with the issue of existence.
We show that
quadratic harnesses are  described by five numerical constants,
which, under appropriate integrability conditions, determine the process
uniquely. We study the related
integrability properties of a slightly wider class of processes, improving
 earlier results of that type
 \cite[Corollary 4]{Bryc85c}, \cite[Theorem 2]{Bryc85b},
\cite[Theorem 2($1^\circ$)]{Wesolowski93}.
We show that martingale polynomials associated with a quadratic harness  with finite moments of all orders
lead to a  $q$-commutation equation
\begin{equation}\label{C-com}
[\xxx,\yyy]_q=\I+\tau \xxx^2+ \sigma\yyy^2+\theta\xxx+ \eta \yyy ,
\end{equation}
where $\xxx,\yyy$ are non-commutative variables,
$$[\xxx,\yyy]_q=\xxx \yyy - q\yyy\xxx,$$
and  $q,\eta,\theta,\sigma,\tau$ are the numerical constants that
describe the quadratic harness. This hints on more connections with non-commutative probability,
$q$-Fock space constructions, and classical versions of non-commutative processes, see
\cite{BKS97}.
In this paper we use \rf{C-com} to  prove that quadratic harnesses with finite moments of all orders
have orthogonal martingale polynomials, and to derive their three step recurrences.
The explicit three terms recurrences that are associated with quadratic harnesses
include a four-parameter family of the 
polynomials in Section \ref{Sect classical}, and a four-parameter family
of $q$-orthogonal polynomials
in Section
\ref{OM sigma-tau=0}.
Based on one of the explicit recurrences from this paper, in \cite{Bryc-Matysiak-Wesolowski-04b} we
 extend the construction from
\cite{Bryc-Wesolowski-04} to the construction of the
quadratic harness which we call
the bi-Poisson processes.

Some of our arguments rely on the symmetries of the problem, but perhaps we did not explore the symmetries deeply enough.
The family of quadratic harnesses that we study is
 invariant under the
 action of  the translations and reflections of $\sR$, and
the associated  affine Hecke algebra is known to be associated with the
  Askey-Wilson polynomials, see \cite{Noumi-Stokman}.
Since the explicit recurrences that we found in Section \ref{SECT: Appl} give
 sub-classes of the reparametrized Askey-Wilson polynomials, it is  plausible that the theory of affine root systems  \cite{Macdonald}
 might lead to additional progress.

The paper is organized as follows.
In Section \ref{SECT: Def} we define quadratic harnesses and state the main results.
Proofs of the main results are given in  Section \ref{SECT: proofs}.
In Section \ref{SECT: Appl} we  derive explicit
recurrences for the orthogonal martingale polynomials. We also introduce an operator
technique that
resembles umbral calculus as the tool that simplifies the proof of the quadratic harness property.


\section{Definitions and main results}\label{SECT: Def}
\subsection{Harnesses and quadratic harnesses} \label{SECT: SQH}
Let $(X_t)_{t> 0}$ be a separable
square
integrable stochastic process with $\sigma$-fields:
 $\Fsu=\sigma\{X_t: t\in(0,s]\cup[u,\infty)\}$.
 Consider the following two functionals
 \begin{eqnarray}
 \Delta_{t,u}= \Delta_{t,u}(X)&=&\frac{X_u-X_t}{u-t} \label{harness a},\\
 \widetilde{\Delta}_{t,u}=\widetilde{\Delta}_{t,u}(X)&=&\frac{uX_t-tX_u}{u-t}. \label{harness b}
\end{eqnarray}
Both expressions appear in the statement of Theorem \ref{T3.1} below and the second one
is the time-inverse  of the first one. Namely, if $(\widetilde{X}_t)_{t>0}$ denotes the time-inverse process
$(t X_{1/t})_{t>0}$, then
$$ \widetilde{\Delta}_{t,u}(X)=\Delta_{1/u,1/t}(\widetilde{X}), \; 0<t<u.
$$

\begin{definition}[\cite{Mansuy-Yor04}] An integrable process $(X_t)_{t> 0}$ is a (simple) harness if
$$\wwo{\Delta_{s,t}}{\Fru}=\Delta_{r,u}$$
for every $r<s<t<u$.
\end{definition}
A trivial example of a  harness is
a Gaussian process with covariance \rf{pre-cov}; additional examples follow Definition \ref{Def QH}.

As pointed out in \cite{Mansuy-Yor04},  the harness  condition is equivalent to the linearity of regression property
\begin{equation}
\label{LH}
\wwo{X_t}{\Fsu}=\aaa X_s+\bbb X_u,
\end{equation}
where the coefficients are given by
\begin{equation}\lbl{EQ: a+b}
\aaa=\frac{u-t}{u-s},\:
\bbb=\frac{t-s}{u-s}.
\end{equation}
It is clear that  $\aaa+\bbb=1$ and $s\aaa  + u\bbb  = t$;
these two identities will often be used.
Since condition \rf{LH} is invariant under the time-inversion $(X_t)\mapsto (t X_{1/t})$,
another equivalent condition for a harness is
$$\wwo{\widetilde{\Delta}_{s,t}}{\Fru}=\widetilde{\Delta}_{r,u}, \, r<s<t<u.$$

The general form of the covariance of a square-integrable harness is as follows.
\begin{proposition}\label{P2.1} If $(X_t)$ is a square-integrable centered harness on $\sR_+$,
then there are constants $c_0, c_1, c_2,c_3$ such that for $s\leq t$ we have
\begin{equation}\label{pre-cov}
\wo(X_tX_s)=c_0+c_1 s+ c_2 t+c_3st.
\end{equation}
\end{proposition}
\begin{proof}
Multiplying \rf{LH} by $X_s$ and averaging we get
$
\wo(X_sX_t)=\aaa \wo(X_s^2)+\bbb \wo(X_sX_u)$. Thus
$\wo(X_sX_t)=\alpha(s)+\beta(s)t$ is linear in $t$ and $\wo(X_tX_u)=\aaa\wo(X_sX_u)+\bbb\wo(X_u^2)=\gamma(u)+\delta(u)t$
is linear in $t$.

This gives a functional equation
$$
\alpha(s)+\beta(s)t=\gamma(t)+\delta(t) s
$$
which, by taking $s=0$, gives $\gamma(t)=\alpha(0)+\beta(0) t$.
Thus
$\wo(X_sX_t)=c_0+c_1s+c_2t+c_3st$.
\end{proof}

For at centered standardized square-integrable $Z$
that is independent of a quadratic harness $(X_t)$, taking $(\tilde{X}_t)=(X_t+(t+1)Z)$ we get a
quadratic harness
with $\wo(\tilde{X}_s\tilde{X}_t)=\wo(X_sX_t)+1+t+s+ts$. To avoid such non-uniqueness,
throughout this paper we assume
\begin{equation}\label{cov}
\wo(X_t)=0,\: \wo(X_tX_s)=\min\{t,s\}.
\end{equation}

 If \rf{cov}
and \rf{LH} hold true  and  $0< s\leq t\leq u$, then
\begin{equation}\label{LR+}
\wwo{X_t}{\calF_{\geq u}}=\frac{t}{u}X_u
\end{equation}
and
\begin{equation}\label{LR-}
\wwo{X_t}{\calF_{\leq s}}=X_s,
\end{equation}
see {\cite[(4) and (5)]{Bryc-Wesolowski-03}}.
From  martingale \rf{LR-} and reverse martingale \rf{LR+} condition it follows that  the limits
$X_0:=\lim_{t\searrow 0} X_t$, and $\lim_{t\to\infty} X_t/t$
exist with probability one. Under assumption \rf{cov} we have \begin{equation}\label{limits}
\lim_{t\searrow 0} X_t=0, \; \lim_{t\to\infty} X_t/t=0,
\end{equation}
 so without loss of generality
we may  extend $(X_t)_{t>0}$ to include the value $X_0=0$ when convenient. Similarly,
we may extend $(\widetilde{X}_t)_{t>0}=(tX_{1/t})_{t>0}$ to include the value $\widetilde{X}_0=0$ corresponding to $t=0$.

We now turn to the quadratic harness condition.
The familiar quadratic martingale condition associated with the martingale property \rf{LR-}
 can be written as
\begin{equation}\label{CV}
\wwo{X_t^2}{{\cal F}_{\leq s}}=X_s^2+t-s.
\end{equation}
This suggests that the quadratic harness condition
associated with the simple harness property \rf{LH} should be written as
\begin{equation}
\label{QH} \wwo{X_t^2}{\Fsu}=Q_{t,s,u}(X_s,X_u),
\end{equation}
where
\begin{equation}\label{Q-form}
Q_{t,s,u}(x,y)=\AAA x^2+\BBB xy+\CCC y^2+\DDD x + \EEE y+\FFF
\end{equation}
 is a quadratic form in variables $x,y$ with time-dependent coefficients; under assumption \rf{cov} we trivially have
\begin{equation}\label{0} s\AAA + s\BBB  +u\CCC +\FFF=t.\end{equation}


\begin{definition}\label{Def QH} A square integrable process $(X_t)_{t> 0}$ is a quadratic harness on $\sR_+$ if it satisfies conditions
  \rf{LH}
and \rf{QH}.
\end{definition}
The well-known examples of quadratic harnesses are the Wiener, Poisson, and Gamma processes.
Ref. \cite{Wesolowski93} identifies all quadratic harnesses with
covariance
\rf{cov}
when the quadratic form on the right hand side of \rf{QH} is such that the corresponding conditional variance
is a function of the increments $X_u-X_s$ only; this adds
to the already listed examples two additional L\'evy processes: the Pascal, and Meixner processes.
A related non-commutative form of condition \rf{QH} appears in \cite{Bozejko-Bryc-04}.
The main result of \cite{Bryc-Wesolowski-03} asserts that quadratic harnesses with covariance \rf{cov}
which satisfy the quadratic martingale condition
\rf{CV} are in fact uniquely determined $q$-Meixner Markov processes.
A  quadratic harness that does not satisfy condition \rf{CV} is analyzed in \cite{Bryc-Wesolowski-04}.

\subsection{Five-parameter representation}
Generically, the quadratic form
on the right hand side of \rf{QH} is determined uniquely up to five numerical parameters.

\begin{theorem}
\label{T3.1}
Let $(X_t)$ be a quadratic harness with covariance \rf{cov}. Suppose that \rf{QH}  holds
with $\Ftsu\ne0$,  for all $0< s<t<u$,
and that
$1,X_s,X_t,X_s X_t, X_s^2,X_t^2$ are linearly independent for all $0< s<t$. Then there exist $\eta,\theta\in\sR$,
 $\sigma,\tau\geq 0$, and $q\leq 1+2\sqrt{\sigma\tau}$  such that
\begin{eqnarray}
&&\wVar{X_t}{\Fsu} =  \nonumber  \\
\label{q-Var}
 &&  \Ftsu\left( 1+ \eta
\widetilde{\Delta}_{s,u} +\theta \Delta_{s,u}+\sigma \widetilde{\Delta}_{s,u}^2
+\tau \Delta_{s,u}^2
-(1-q)\Delta_{s,u}\widetilde{\Delta}_{s,u} \right)
\end{eqnarray}
for all $0< s<t<u$, where
\begin{equation}
\Ftsu = \frac{(u-t)(t-s)}{u(1+\sigma s)+\tau-q s}\label{F0}\end{equation}
is the normalizing function, and $\Delta_{s,u}$, $\widetilde{\Delta}_{s,u}$ are given by (\ref{harness a}-\ref{harness b}).
\end{theorem}
Recall that
the conditional variance of $X$  with respect to a $\sigma$-field
$\calF$ is defined as
$\wVar{X}{\calF}=\wwo{X^2}{\calF}-\left(\wwo{X}{\calF}\right)^2$.

\begin{remark}
\label{Time inversion}
Time-inversion $(X_t)\mapsto(tX_{1/t})$
preserves the class of quadratic harnesses, modifying the coefficients in 
\rf{QH}.
More precisely, suppose $(X_t)$ satisfies the assumptions of Theorem \ref{T3.1}, and let $\widetilde{X}_t=tX_{1/t}$ be its time inverse.
 Then
 $(\widetilde{X}_t)$  is a quadratic harness with respect to its $\sigma$-fields
$\widetilde{\mathcal F}_{s,u}={\mathcal F}_{1/u,1/s}$, and \rf{q-Var} holds with the roles of the parameters
 $(\eta,\theta)$ and $(\sigma,\tau)$ switched within each pair:
\[
\wVar{\widetilde{X}_t}{\widetilde{\mathcal F}_{s,u}}/ \widetilde{F}_{t,s,u}\]
\[=
1+ \theta
\widetilde{\Delta}_{s,u}(\widetilde{X}) +\eta \Delta_{s,u}(\widetilde{X})+\tau \widetilde{\Delta}_{s,u}^2(\widetilde{X})
+\sigma \Delta_{s,u}^2(\widetilde{X})
-(1-q)\Delta_{s,u}(\widetilde{X})\widetilde{\Delta}_{s,u}(\widetilde{X}),
\]
where $\widetilde{F}_{t,s,u}=\frac{(u-t)(t-s)}{u(1+\tau
s)+\sigma-qs}$.


Ref. \cite{Gallardo-Yor04b}  gives criteria for time-inversion invariance of Markov processes.
\end{remark}

\subsection{Orthogonal martingale polynomials}
Suppose  that a quadratic harness $(X_t)$ has moments of all orders and martingale polynomials $p_n(x;t)$
of all degrees $n\geq 0$, that is
\begin{equation}\label{ompind}
\wwo{\pln{n}{t}}{\mathcal{F}_{\le s}}=\pln{n}{s}, \; 0<s<t.
\end{equation}
Clearly,  $p_0=1$ and $p_1(x;t)=x$ are natural initial choices, see  \rf{LR-}.
Since $xp_n(x;t)$ is a polynomial of degree $n+1$, it follows that
\begin{equation}\label{p-rec}
x p_n(x;t)=\sum_{k=0}^{n+1} C_{k,n}(t) p_k(x;t).
\end{equation}


\begin{theorem}
\label{T5.1}  Suppose a quadratic harness  $(X_t)$
with covariance \rf{cov} and
conditional variance \rf{q-Var}
has finite moments of all orders
and martingale polynomials $p_n(x;t)$.
If for each $t>0$ the random variable $X_t$ has infinite
support, then
   recurrence \rf{p-rec}  holds with the infinite matrices
\[
\mC_t:=\left[ {\begin{array}{*{30}l}
   C_{00}(t) & C_{01}(t) & C_{02}(t) & C_{0,3}(t) & \dots\vspace{.25cm} &  \\
   C_{10}(t) & C_{11}(t) & C_{12}(t) & C_{13}(t) & \dots\vspace{.25cm} & \\
   0 & C_{21}(t) & C_{22}(t) & C_{23}(t) & \dots\vspace{.25cm}  & \\
   0 & 0 & C_{32}(t) & C_{33}(t) & \dots\vspace{.1cm} & \\
   0 &  0 &  0 &  C_{43}(t) & \ddots\vspace{.0cm}  & \\
    \vdots &  \vdots &  \vdots &   & {} &\ddots \\
 \end{array} } \right]
\]
   given by
\begin{equation}\label{C-lin}
\mC_t=t\xxx+\yyy, \;t> 0,
\end{equation} %
and the infinite matrices $\xxx=\mC_1-\mC_0$, $\yyy=\mC_0$ satisfy  equation \rf{C-com}.
\end{theorem}


From Theorem  \ref{T5.1} we derive a number of  equations that eventually determine
 the orthogonal martingale polynomials. Namely, it is well known that
 for orthogonal martingale polynomials $p_n(x;t)$
  recurrence \rf{p-rec} holds with a  tri-diagonal matrix $\mC$, see \cite{Chihara}.
Writing \rf{p-rec} as
\begin{equation}\lbl{three step abstract}
x p_n(x;t)=a_{n}(t)p_{n+1}(x;t)+b_n(t)p_n(x;t)+c_n(t)p_{n-1}(x;t), \; n\geq 0,
\end{equation}
from \rf{C-lin} we get
\begin{equation}\label{coeff}
a_n(t)=\sigma\alpha_{n+1} t+\beta_{n+1}\;,\;\;\;b_n(t)=\gamma_n t+\delta_n\;,\;\;\;c_n(t)=\eps_n t+\varphi_n,
\end{equation}
and  the coefficients in \rf{coeff}
satisfy a number of equations that result from the $q$-commutation equation \rf{C-com}.

Setting  $p_{-1}(x;t)=0$, $p_0(x;t)=1$, $p_1(x;t)=x$, and using \rf{cov} we see that
the initial values are given by
\begin{equation}\label{ini}
\alpha_1=0,\; \beta_1=1,\;\gamma_0=\delta_0=0,\;
\eps_1=1,\, \varphi_1=0.
\end{equation}
For $n\geq 1$, \rf{C-com} implies that coefficients $\alpha_n,\beta_n,\gamma_n,\delta_n,\varphi_n,\eps_n$ satisfy
  \begin{equation}\label{eq1}
\sigma^2\tau\alpha_n\alpha_{n+1} +\sigma\alpha_n\beta_{n+1}q+\sigma\beta_n\beta_{n+1}=\sigma\alpha_{n+1}\beta_n,
\end{equation}
\begin{eqnarray}\label{eq2}
&&\beta_{n+1}\gamma_{n+1}+\sigma\alpha_{n+1}\delta_n
\\ \nonumber
&&=\sigma\alpha_{n+1}(\gamma_n+\gamma_{n+1})\tau+(\sigma\alpha_{n+1}\delta_{n+1}+\beta_{n+1}\gamma_n)q
+\beta_{n+1}(\delta_n+\delta_{n+1})\sigma \\ \nonumber &&+\sigma\alpha_{n+1}\theta+\beta_{n+1}\eta,
\end{eqnarray}
\begin{eqnarray}\label{eq3}
&&\beta_{n+1}\varepsilon_{n+1}+\gamma_n\delta_n+\sigma\alpha_{n}\varphi_n
\\ \nonumber
&&=(\sigma\alpha_{n+1}\varepsilon_{n+1}+\gamma_n^2+\sigma\alpha_{n}\varepsilon_n)\tau
+(\sigma\alpha_{n+1}\varphi_{n+1}+\gamma_n\delta_n+\beta_{n}\varepsilon_n)q
\\ \nonumber
&&+(\beta_{n+1}\varphi_{n+1}+\delta_n^2+\beta_{n}\varphi_n)\sigma+\gamma_n\theta+\delta_n\eta+1,
\end{eqnarray}

\begin{eqnarray}\label{eq4}
&&\gamma_{n-1}\varphi_n+\delta_n\varepsilon_n \\ \nonumber
&&=(\gamma_{n-1}+\gamma_n)\varepsilon_n\tau+(\gamma_n\varphi_n+\delta_{n-1}\varepsilon_n)q
+(\delta_{n-1}+\delta_n)\varphi_n\sigma+\varepsilon_n\theta+\varphi_n\eta,
\end{eqnarray}
\begin{equation}\label{eq5}
\varepsilon_{n}\varphi_{n+1}=\varepsilon_{n}\varepsilon_{n+1}\tau+\varepsilon_{n+1}\varphi_{n}q+\varphi_{n}\varphi_{n+1}\sigma.
\end{equation}

If $a_n(t)\ne 0$ for all $n$, then %
 the three step recurrence \rf{three step abstract}
defines a family of polynomials in variable $x$.
We now show that these are indeed the
martingale orthogonal polynomials  for $(X_t)$.

\begin{theorem}\label{C mo poly}
Suppose a quadratic harness  $(X_t)$
with covariance \rf{cov} and conditional variance \rf{q-Var}
has finite moments of all orders, parameters $q, \eta,\theta,\sigma,\tau$ are such that
$(q,\sigma\tau)\ne (-1,1)$, and
 the equations (\ref{ini}-\ref{eq5})
have a solution such that $a_n(t) \ne 0$  for all $t\geq 0$, $n\geq 0$.
Let  $\{p_n(x;t):n\geq 0, t>0\}$
satisfy 
\rf{three step abstract} and \rf{coeff}  for $n=0,1,...$
with  $p_{-1}(x,t)=0, p_0(x,t)=1$, $p_1(x,t)=x$. Then
 $\{p_n(x;t)\}$ are  orthogonal martingale
 polynomials for $(X_t)$.
\end{theorem}

 In Section \ref{SECT: Appl}, we give sufficient conditions in terms of parameters $q, \eta,\theta,\sigma,\tau$ for
  the assumptions of Theorem \ref{C mo poly} to be satisfied, and
  we derive explicit three step recurrences for the orthogonal martingale polynomials in several cases of interest.

Theorem 
\ref{C mo poly}
implies that if   quadratic harnesses $(X_t)_{t>0}$ and $(Y_t)_{t>0}$ satisfy its assumptions
with the same parameters $\eta,\theta,\sigma,\tau,q$, then
\[
\wo(X_{t_1}^{n_1}X_{t_2}^{n_2}\dots X_{t_k}^{n_k})=\wo(Y_{t_1}^{n_1}Y_{t_2}^{n_2}\dots Y_{t_k}^{n_k})
\]
 for all  $k\geq 1$,
$0<t_1<t_2<\dots<t_k$, and $n_1,\dots,n_k\in\sN$.
Under appropriate integrability conditions, this implies that $(X_t)$ is Markov and is
 uniquely determined by the  parameters $\eta,\theta,\sigma,\tau,q$. In particular, this is the case when
 $\sup_n|b_n(t)|<\infty$ and $\sup_n |a_{n-1}(t)c_n(t)|<\infty$, as in this case
the recurrence \rf{three step abstract} corresponds to a compactly supported measure, see
 \cite[Section 2]{Askey-Ismail84MAMS}.
 However,  the question of  existence of such a process is non-trivial and the constructions
  are known only in special cases, see
\cite{Bryc-Matysiak-Wesolowski-04b}, \cite{Bryc-Wesolowski-04}, \cite{Bryc-Wesolowski-03}.

Next, we show that the integrability assumption of Theorem 
\ref{C mo poly}
is automatically
 satisfied if $\sigma\tau=0$.
We remark that the result stated below does not use the full power of the quadratic harness condition,
 and generalizes
 \cite[Theorem 2]{Bryc85b} and \cite[Theorem 2($1^\circ$)]{Wesolowski93}.
Since
$$\wVar{X_t}{\Fsu}=\wwo{X_t^2}{\Fsu}-(\atsu X_s+\btsu X_u)^2,$$
using \rf{limits} we can pass to the limit  in \rf{q-Var}
as $u\to\infty$ or as $s\to0$. This gives
\begin{eqnarray} \label{CV-}
\wVar{X_t}{\calF_{\leq s}}&=&\frac{t-s}{1+\sigma s}\left(\sigma X_s^2+
\eta X_s+1\right),\\ \label{RV-}
\wVar{X_t}{\calF_{\geq u}}&=&\frac{t(u-t)}{u+\tau}\left(\tau \frac{X_u^2}{u^2}+\theta \frac{X_u}{u}+1\right).
\end{eqnarray}

\begin{theorem}\label{T4.1} 
If  a square-integrable stochastic process $(X_t)$ with covariance \rf{cov}
 satisfies  \rf{LR+}, \rf{LR-},  \rf{CV-} and \rf{RV-} with $\sigma,\tau\geq 0$ such that
$ \sigma\tau \leq  1/2^{4r+10}$ for some $r> 2$
 then $\wo(|X_t|^{r})<\infty$ for all $t>0$. In particular, if $\sigma\tau=0$, then
 $\wo(|X_t|^{r})<\infty$ for all $r,t>0$.
\end{theorem}

\section{Proofs of the main results}\label{SECT: proofs}

\subsection{Proof of Theorem \ref{T3.1}}
For $0< x<y<z$ define 
\begin{eqnarray*}
\sigma_{y,x,z}:=&&\frac{A_{y,x,z}+B_{y,x,z}+C_{y,x,z}-1}{F_{y,x,z}}\ ,\\
\tau_{y,x,z}:=&&\frac{x^2 A_{y,x,z}+x z B_{y,x,z}+z^2 C_{y,x,z}-y^2}{F_{y,x,z}}\ ,\\
1+q_{y,x,z}:=&&\frac{B_{y,x,z}(z-x)}{F_{y,x,z}}\ ,\\
\eta_{y,x,z}:=&&\frac{D_{y,x,z}+E_{y,x,z}}{F_{y,x,z}}\ ,\\
\theta_{y,x,z}:=&&\frac{x D_{y,x,z}+z E_{y,x,z}}{F_{y,x,z}}\ .
\end{eqnarray*}

We will show that the left hand sides of the above equations do not depend on the arguments
$y,x,z$. To this end we prove three claims.
Throughout the proofs of all three Claims, 
$0< r<s<t<u$ are arbitrary numbers.
\begin{claim}
\label{claim 1} For all $0< x<y<z$, and for $f=\sigma,\tau,q,\eta,\theta$ we have
\begin{equation}\label{nogasrodkowa}
f_{y,x,z}=f_{\widetilde{y}, x, z},
\end{equation}
provided
$x<\widetilde{y}<z$.
\end{claim}
\begin{proof}[Proof of Claim \ref{claim 1}]

Observe that
from (\ref{LH}) and (\ref{QH}) we get
\[
\wwo{X_s X_t}{\Fru}=\wwo{X_s \wwo{X_t}{\Fsu}}{\Fru}
\]
\[=\atsu\wwo{X_s^2}{\Fru}
+\btsu X_u\wwo{X_s}{\Fru}
\]
\[=\atsu\left(\Asru X_r^2+\Bsru X_rX_u+\Csru X_u^2+\Dsru X_r+
\Esru X_u+\Fsru\right)
\]
\[+\btsu X_u\left(\asru X_r+\bsru X_u\right)=
\atsu \Asru X_r^2+\left(\atsu \Bsru+\btsu\asru\right)X_rX_u
\]
\[+
\left(\atsu \Csru+\btsu\bsru \right)X_u^2
+\atsu \Dsru X_r+\atsu \Esru X_u+\atsu \Fsru\;.
\]

On the other hand,
\[
\wwo{X_s X_t}{\Fru}=\wwo{\wwo{X_s}{\Frt}X_t}{\Fru}
=\asrt X_r\wwo{X_t}{\Fru}+\bsrt \wwo{X_t^2}{\Fru}\]
\[
=\asrt X_r\left(\atru X_r+\btru X_u\right)\]
\[+
\bsrt \left(\Atru X_r^2+\Btru X_rX_u+\Ctru X_u^2\right.
\left.+\Dtru X_r+
\Etru X_u+\Ftru\right)
\]
\[
=\left(\asrt \atru +\bsrt \Atru\right)X_r^2+
\left(\asrt\btru+\bsrt \Btru\right)X_rX_u\]
\[
+\bsrt \Ctru X_u^2+
\bsrt \Dtru X_r+\bsrt \Etru X_u+\bsrt \Ftru\;.
\]

Comparing the coefficients 
at $X_r^2$, $X_rX_u$, $X_u^2$, $X_r$,
$X_u$ and 1 in the above expressions we get
\begin{eqnarray}
\atsu\Asru &=& \bsrt\Atru+\asrt\atru\ , \label{A1}\\
\atsu\Bsru &=& \bsrt\Btru\ , \label{B1}\\
\btsu\bsru+\atsu\Csru &=& \bsrt\Ctru\ , \label{C1}\\
\atsu\Dsru &=& \bsrt\Dtru\ , \label{D1}\\
\atsu\Esru &=& \bsrt\Etru\ , \label{E1}\\
\atsu\Fsru &=& \bsrt\Ftru\ .  \label{F1}
\end{eqnarray}
(
In the derivation of (\ref{B1}) we used the fact that $\btsu\asru=\asrt\btru$).

Adding (\ref{A1}), (\ref{B1}) and (\ref{C1}),
we get
\[
\atsu\left(\Asru+\Bsru+\Csru-1\right)+\atsu+\btsu\bsru\]
\[=\bsrt\left(\Atru+\Btru+\Ctru-1\right)+\bsrt+\asrt\atru.
\]
Since a calculation shows that \begin{equation}\label{tmp}
\atsu+\btsu\bsru=\bsrt+\asrt\atru,\end{equation}  dividing by (\ref{F1}), we get
$\sigma_{s,r,u}=\sigma_{t,r,u}$ which proves \rf{nogasrodkowa} when $f=\sigma$.

We proceed similarly when $f=\tau,q,\eta,\theta$. Adding (\ref{A1}) multiplied by $r^2$, (\ref{B1}) multiplied by
$ru$ 
(\ref{C1}) multiplied by $u^2$, and dividing by (\ref{F1}), we obtain
$\tau_{s,r,u}=\tau_{t,r,u}$, 
after
noticing that
$
s^2\atsu + u^2\btsu\bsru =
  t^2\bsrt + r^2\asrt\atru
$. Equation (\ref{B1}) multiplied by $(u-r)$ and divided by (\ref{F1})
gives $q_{s,r,u}=q_{t,r,u}$. Adding equations (\ref{D1}) and (\ref{E1}) and dividing by (\ref{F1})
gives $\eta_{s,r,u}=\eta_{t,r,u}$ and similarly  multiplying (\ref{D1}) by $r$ and (\ref{E1}) by
$u$ after dividing by (\ref{F1}) gives $\theta_{s,r,u}=\theta_{t,r,u}$. Thus we obtained
$f_{s,r,u}=f_{t,r,u}$ for $f=\sigma,\tau,q,\eta,\theta$.
If $0< x<y<z$ then substitution $r=x,s=y,u=z$ yields (\ref{nogasrodkowa}) for
$\widetilde{y}\in(y,z)$; substitution $r=x,t=y,u=z$ gives (\ref{nogasrodkowa}) for
$\widetilde{y}\in(x,y)$, completing the proof of Claim (\ref{claim 1}).
\end{proof}

\begin{claim}\label{claim 2}
For all $0<x<y<z$, and $f=\sigma,\tau,q,\eta,\theta$ we have
\begin{equation}\label{nogalewa}
f_{y,x,z}=f_{y,\widetilde{x},z},
\end{equation}
provided $0<\widetilde{x}<y$.
\end{claim}
\begin{proof}[Proof of Claim \ref{claim 2}]
Now consider the identity
\[
\wwo{X_t^2}{\Fru}=\wwo{\wwo{X_t^2}{\Fsu}}{\Fru}.
\]
Using (\ref{LH}) and (\ref{QH}) we see that $\wwo{\wwo{X_t^2}{\Fsu}}{\Fru}$ is given by
\[
\wwo{\Atsu X_s^2+\Btsu X_sX_u+\Ctsu X_u^2+\Dtsu X_s+\Etsu X_u+\Ftsu}{\Fru}
\]
\[=\Atsu\wwo{X_s^2}{\Fru}
+\Btsu X_u\wwo{X_s}{\Fru}+\Ctsu X_u^2+\Dtsu
\wwo{X_s}{\Fru}
\]
\[
+\Etsu X_u+\Ftsu=\Atsu\big(\Asru X_r^2+\Bsru X_rX_u+\Csru X_u^2+\Dsru X_r\]
\[+
\Esru X_u+\Fsru\big)+\Btsu X_u\left(\asru X_r+\bsru X_u\right)+\Ctsu X_u^2
\]
\[+\Dtsu \left(\asru X_r+\bsru X_u\right)+\Etsu X_u+
\Ftsu=
\Atsu \Asru X_r^2\]
\[+\left(\Atsu\Bsru+\Btsu\asru \right)X_rX_u+
\left(\Atsu\Csru+\Btsu\bsru +\Ctsu\right)X_u^2\]
\[+
\left(\Atsu\Dsru+\Dtsu\asru \right)X_r+(\Atsu\Esru+\Dtsu\bsru+\Etsu)X_u\]
\[+
\Atsu\Fsru+\Ftsu.
\]
On the other hand,
\[
\wwo{X_t^2}{\Fru}=\Atru X_r^2+\Btru X_rX_u+\Ctru X_u^2+\Ctru X_r+\Dtru X_u+\Ftru.
\]

Comparing 
the  coefficients at $X_r^2$, $X_r X_u$, $X_u^2$, $X_r$,
$X_u$ and $1$, we obtain
\begin{eqnarray}
\Atru &=& \Atsu\Asru \label{A2}\ ,\\
\Btru &=& \Atsu\Bsru+\Btsu\asru\ , \label{B2}\\
\Ctru &=& \Atsu\Csru+\Btsu\bsru+\Ctsu\ , \label{C2}\\
\Dtru &=& \Atsu\Dsru+\Dtsu\asru\ , \label{D2}\\
\Etru &=& \Atsu\Esru+\Dtsu\bsru+\Etsu\ , \label{E2}\\
\Ftru &=& \Atsu\Fsru+\Ftsu\ . \label{F2}
\end{eqnarray}
Substituting the right hand sides of equation (\ref{A2}-\ref{F2}) for
$\Atru,\Btru,\dots,\Ftru$ on the right hand sides of  (\ref{A1}-\ref{F1}) we get
\begin{align}
\Asru(\atsu-\bsrt\Atsu) &=\asrt\atru\ , &\label{A3}\\
\Bsru(\atsu-\bsrt\Atsu) &= \bsrt\Btsu\asru\ , &\label{B3}\\
\Csru(\atsu-\bsrt\Atsu) &= \bsrt(\Btsu\bsru+\Ctsu)-\btsu\bsru\ ,& \label{C3}\\
\Dsru(\atsu-\bsrt\Atsu) &= \bsrt\Dtsu\asru\ , &\label{D3}\\
\Esru(\atsu-\bsrt\Atsu) &= \bsrt(\Dtsu\bsru+\Etsu)\ , &\label{E3}\\
\Fsru(\atsu-\bsrt\Atsu) &= \bsrt\Ftsu\ . &\label{F3}
\end{align}
We can now proceed analogously to the proof of Claim \ref{claim 1}. Namely, adding (\ref{A3}),
(\ref{B3}), (\ref{C3}), again taking into account
\rf{tmp}, and dividing by (\ref{F3}), we get $\sigma_{s,r,u}=\sigma_{t,s,u}$.
Adding (\ref{A3}) multiplied by $r^2$, (\ref{B3}) multiplied by $ru$ and (\ref{C3}) multiplied by
$u^2$, and dividing by (\ref{F3}), we obtain $\tau_{s,r,u}=\tau_{t,s,u}$. Equation (\ref{B3})
multiplied by $(u-r)$ and divided by (\ref{F3}) gives $q_{s,r,u}=q_{t,s,u}$. Adding equations
(\ref{D3}) and (\ref{E3}) and dividing by (\ref{F3}) gives $\eta_{s,r,u}=\eta_{t,s,u}$ and
similarly  multiplying (\ref{D3}) by $r$ and (\ref{E3}) by $u$ after dividing by (\ref{F3}) gives
$\theta_{s,r,u}=\theta_{t,s,u}$. Thus we obtained $f_{s,r,u}=f_{t,s,u}$ for
$f=\sigma,\tau,q,\eta,\theta$. By Claim \ref{claim 1}
$f_{s,r,u}=f_{t,r,u}$, so  $f_{t,r,u}=f_{t,s,u}$,
ending the proof of
Claim \ref{claim 2}.
\end{proof}

\begin{claim}\label{claim 3}
For all $0< x<y<z$, and $f=\sigma,\tau,q,\eta,\theta$ we have
\begin{equation}\label{nogaprawa} 
f_{y,x,z}=f_{y,x,\widetilde{z}},
\end{equation}
provided $\widetilde{z}>y$.
\end{claim}
\begin{proof}[Proof of Claim \ref{claim 3}]
This follows from Claim \ref{claim 2} by the time-inversion $(X_t)\mapsto (t X_{1/t})$.
Alternatively, one can repeat the previous arguments, starting with
 the identity
\[
\wwo{X_s^2}{\Fru}=\wwo{\wwo{X_s^2}{\Frt}}{\Fru}.
\]
\end{proof}

\begin{proof}[Conclusion of proof of Theorem \ref{T3.1}]
Now, it is easy to deduce that
(\ref{nogasrodkowa}), (\ref{nogalewa}) and (\ref{nogaprawa}) imply that functions
$\sigma,\tau,q,\eta,\theta$ are in fact constants. Indeed, given
$0< x_1<y_1<z_1$ and $0< x_2<y_2<z_2$, and $f=\sigma,\tau,q,\eta,\theta$, with
$\widetilde{x}:=
\min\left\{x_1,x_2\right\}$ and $\widetilde{z}:=\max\left\{z_1,z_2\right\}$, we see
that
\begin{align*}
f_{y_1,x_1,z_1}& =f_{y_1,x_1,\widetilde{z}} && \text{by Claim \ref{claim 3},}\\
& =f_{y_1,\widetilde{x},\widetilde{z}} && \text{by Claim \ref{claim 2},}\\
& =f_{y_2,\widetilde{x},\widetilde{z}} && \text {by Claim \ref{claim 1},}\\
& =f_{y_2,x_2,\widetilde{z}} && \text {by Claim \ref{claim 2},}\\
& =f_{y_2,x_2,z_2} && \text {by Claim \ref{claim 3}.}\\
\end{align*}
(If $x_2=\widetilde{x}$ or $z_2=\widetilde{z}$, one or both of the last two steps is unnecessary.)

So, for all $0< s<t<u$ we have
\begin{gather}
\sigma=\frac{\Atsu+\Btsu+\Ctsu-1}{\Ftsu}\ ,
\tau=\frac{s^2 \Atsu+s u \Btsu+u^2 \Ctsu-t^2}{\Ftsu}\ , \label{sigma-tau}\\
1+q=\frac{\Btsu(u-s)}{\Ftsu}\ , \eta=\frac{\Dtsu+\Etsu}{\Ftsu}\ , \theta=\frac{s
\Dtsu+u\Etsu}{\Ftsu}\ . \label{q-eta-theta}
\end{gather}
The above equations, along with (\ref{0}), make a system of linear equations in variables
$\Atsu,\ldots,\Ftsu$. This system must be solvable for all $0<s<t<u$, so its determinant
$(u-s)^3(u(1+\sigma s)+\tau-qs)\ne 0$.
As $u\to\infty, s\to 0$, the expression $u(1+\sigma s)+\tau-qs$ is positive. Thus
$u(1+\sigma s)+\tau-qs>0$. Taking  the limits $s\to 0, u\to 0$, $s=1,u\to\infty$, and $s\to u$, we get
$\tau\geq 0$, $\sigma\geq 0$ and $q\leq 1+\sigma u + \tau/u$ respectively. Minimizing the latter over $u>0$
we get $q\leq 1+2\sqrt{\sigma\tau}$.

The unique solution of the system of equations is given by \rf{F0},  (\ref{A0}-\ref{E0}). A
calculation shows that (\ref{q-Var}) holds.
\end{proof}

\begin{remark}
A calculation shows that formula \rf{q-Var} is equivalent to the following formulas, valid for all $0<s<t<u$

\begin{eqnarray}
\Atsu &=& \frac{(u-t)[u(1+\sigma t)+\tau-q t]}{(u-s)[u(1+\sigma
s)+\tau-q s]}\ ,\label{A0}\\
\Btsu &=& \frac{(u-t)(t-s)(1+q)}{(u-s)[u(1+\sigma s)+\tau-q s]}\ ,\label{B0}\\
\Ctsu &=& \frac{(t-s)[t(1+\sigma s)+\tau-q s]}{(u-s)[u(1+\sigma s)+\tau-q s]}\ ,\label{C0}\\
\Dtsu &=& \frac{(u-t)(t-s)(u\eta-\theta)}{(u-s)[u(1+\sigma
s)+\tau-q s]}\ ,\label{D0}\\
\Etsu &=& \frac{(u-t)(t-s)(\theta-s\eta)}{(u-s)[u(1+\sigma
s)+\tau-q s]}\ ,\label{E0}
\end{eqnarray}
and $\Ftsu $ is given by \rf{F0}.
\end{remark}
\subsection{Proof of Theorem \ref{T5.1}}
For non-commutative variables $\xxx,\yyy$, let
\begin{equation}\label{NC Q def}
Q_{t,s,u}(\xxx,\yyy)=\Atsu \xxx^2+\Btsu \xxx\yyy +\Ctsu \yyy^2+\Dtsu \xxx + \Etsu \yyy +\Ftsu;
\end{equation} recall
 the quadratic form \rf{Q-form}  in commuting variables $x,y$.
(Note that the order of $\xxx,\yyy$ at $\Btsu$ is changed when compared to its dual version \rf{NC Q*}.)

 \begin{lemma}\label{Lq 1} Under the assumptions of Theorem \ref{T5.1},
 \begin{eqnarray}
 \label{NC L}
 \mC_t=\atsu\mC_s+\btsu\mC_u\,, \\
 \label{NC Q}
\mC_t^2=Q_{t,s,u}(\mC_s,\mC_u)\,,
 \end{eqnarray}
 where the coefficients in \rf{NC Q def} are given by \rf{F0} and (\ref{A0}-\ref{E0}).
 \end{lemma}
\begin{proof}

For a  polynomial $\varphi:\sR\to\sR$ consider vectors
$$\v_{s,t}=[\wo(\varphi(X_s)p_0(X_t;t)), \wo(\varphi(X_s)p_1(X_t;t)),\dots,\wo(\varphi(X_s)p_k(X_t;t)),\dots]\in\sR^\infty.$$
Since for $s>0$ random  variable $X_s$ has infinite support, polynomials $1, X_s, X_s^2, \dots$ are linearly independent,
and the corresponding orthogonal polynomials are non-degenerate.
This implies that as we change $\varphi$ for a fixed $s>0$,
vectors of the form $\v_{s,s}$ are dense in $\sR^\infty$, equipped with the product topology.
Indeed, applying the Gram-Schmidt orthogonalization process to
$p_0(x;s),p_1(x;s),\dots$ we get a sequence of orthogonal polynomials $q_0,q_1,\dots$ such that $\wo(q_k(X_s)^2)\ne 0$.
Therefore, for any $\w=[w_0,w_1,\dots]\in\sR^\infty$, and $n\geq 0$ we can find numbers $u_0,u_1,\dots,u_n$ such that
 the first $n$ coordinates of $\v_{s,s}$ corresponding to $\varphi=\sum_{j=0}^n u_jq_j$  are equal to $w_0,w_1,\dots,w_n$.

This implies that to verify the identities, it suffices to verify that the identities hold true when multiplied from the left by $\v_{s,s}$.

By the  martingale polynomial property, $\wo(\varphi(X_s)p_k(X_s;s))=\wo(\varphi(X_s)p_k(X_t;t))$. Therefore, from
\rf{p-rec} it follows that $\v_{s,s}=\v_{s,t}$ and
\[
\v_{s,s}\times \mC_t=\v_{s,t}\times \mC_t\]
\[
=\left[\wo\left(\varphi(X_s)\sum_jC_{j,0}(t)p_j(X_t;t)\right),\wo\left(\varphi(X_s)\sum_jC_{j,1}(t)p_j(X_t;t)\right),\dots  \right]
\]
\[=\left[\wo(\varphi(X_s)X_tp_0(X_t;t)), \wo(\varphi(X_s)X_tp_1(X_t;t)),\dots,\wo(\varphi(X_s)X_tp_k(X_t;t)),\dots\right]
\]
\[=\left[\wo(\varphi(X_s)X_tp_0(X_u;u)), \wo(\varphi(X_s)X_tp_1(X_u;u)),\dots,\wo(\varphi(X_s)X_tp_k(X_u;u)),\dots\right]
\]
\[=\atsu \left[\wo(\varphi(X_s)X_sp_k(X_u;u)): k\geq0\right]
+\btsu \left[\wo(\varphi(X_s)X_up_k(X_u;u)): k\geq0\right].
\]
Using the  martingale polynomial property again, we see that $\v_{s,s}\times \mC_t$ equals to
\[\atsu \left[\wo(\varphi(X_s)X_sp_k(X_s;s)): k\geq0\right]+\btsu \left[\wo(\varphi(X_s)X_up_k(X_u;u)): k\geq0\right]\]
\[\atsu \v_{s,s}\times\mC_s+\btsu\v_{s,u}\times \mC_u=\v_{s,s}\times(\atsu \mC_s+\btsu \mC_u),
\]
proving \rf{NC L}. Similar reasoning proves \rf{NC Q}:
\[
\v_{s,s}\times \mC_t^2=\v_{s,t}\times \mC_t^2=
\left[\wo(\varphi(X_s)X_tp_k(X_t;t)): k\geq0\right]\times \mC_t\]
\[=
\left[\wo(\varphi(X_s)X_t^2p_k(X_t;t)):  k\geq0\right]
=\left[\wo(\varphi(X_s)X_t^2p_k(X_u;u)): k\geq0\right]\]
\[
=\left[\wo(\varphi(X_s)\wwo{X_t^2}{\Fsu}p_k(X_u;u)): k\geq0\right]\]
\[
=\left[\wo(\varphi(X_s)Q_{t,s,u}(X_s,X_t)p_k(X_u;u)): k\geq0\right]
\]
\[
=\Atsu\left[\wo(\varphi(X_s)X_s^2p_k(X_s;s)): k\geq0\right]+
\Btsu\left[\wo(\varphi(X_s)X_sX_up_k(X_u;u)): k\geq 0\right]
\]
\[+\Ctsu \left[\wo(\varphi(X_s)X_u^2p_k(X_u;u)):k\geq 0\right]
+ \Dtsu \left[\wo(\varphi(X_s)X_s p_k(X_s;s)):k\geq 0\right]
\]
\[+\Etsu \left[\wo(\varphi(X_s)X_u p_k(X_u;u)):k\geq 0\right]
+\Ftsu \left[\wo(\varphi(X_s) p_k(X_s;s)):k\geq 0\right]
\]
\[
=\Atsu \v_{s,s}\times \mC_s^2+\Btsu \left[\wo(\varphi(X_s)X_s p_k(X_s;s)):k\geq 0\right]\times \mC_u
\]
\[
+ \Ctsu \v_{s,s}\times\mC_u^2+\Dtsu\v_{s,s}\times\mC_s+\Etsu\v_{s,s}\times\mC_u+\Ftsu \v_{s,s}
\]
\[
=\v_{s,s}\times \left(
\Atsu  \mC_s^2+\Btsu \mC_s\mC_u
+ \Ctsu \mC_u^2+\Dtsu \mC_s+\Etsu \mC_u+\Ftsu \II\right).
\]
\end{proof}
\begin{lemma}\label{Lq 2} Suppose a linear operator $\mC(t)=t\xxx+\yyy$ satisfies relation \rf{NC Q} for some quadratic form
$Q_{t,s,u}$ given by \rf{NC Q def}.
If the coefficients of $Q_{t,s,u}$ are given by \rf{F0}, (\ref{A0}-\ref{E0}) then
\rf{C-com} holds true.
\end{lemma}
\begin{proof}
Expanding the expressions like $\mC_t^2=t^2 \xxx^2+t(\xxx\yyy+\yyy\xxx)+\yyy^2$ and using \rf{NC Q} we get
\begin{eqnarray}\label{pre-q-rel}
t^2 \xxx^2+t(\xxx\yyy+\yyy\xxx)+\yyy^2=
\Atsu (s^2 \xxx^2+s(\xxx\yyy+\yyy\xxx)+\yyy^2)
\\ \nonumber+
\Btsu(\yyy^2+u\yyy\xxx+s\xxx\yyy+su\xxx^2)
+\Ctsu(u^2 \xxx^2+u(\xxx\yyy+\yyy\xxx)+\yyy^2)
\\ \nonumber
+(\Dtsu+\Etsu)\yyy+(s\Dtsu+u\Etsu)\xxx+\Ftsu\I.
\end{eqnarray}
Applying the relations \rf{0}, \rf{sigma-tau}, and \rf{q-eta-theta} to the coefficients at each of the monomials in \rf{pre-q-rel}
we get
\[
\xxx\yyy+\yyy\xxx=\tau\xxx^2+\sigma\yyy^2+\eta\yyy+\theta\xxx+\I+(1+q)\yyy\xxx,
\]
which is the same as \rf{C-com}.
 \end{proof}

\subsection{Proof of Theorem \ref{C mo poly} 
}
First we shall show  the martingale property \rf{ompind}.
To this end, we proceed by induction on $n$. Trivially, \eqref{ompind} holds for
$n=0,1$, see \rf{LR-}.
Assume then that  $n>1$ and \eqref{ompind} holds 
for $0, 1,\ldots,n-1$
and all $0<s<t$.

We start from calculating $\wwo{X_t \pln{n-1}{u}}{\mathcal{F}_{\le s}}$ in two ways. On one hand

\begin{eqnarray*}
& &\wwo{X_t \pln{n-1}{u}}{\mathcal{F}_{\le
s}}=\wwo{\wwo{X_t}{\mathcal{F}_{s,u}}\pln{n-1}{u}}{\mathcal{F}_{\le s}}=\\
&=&\wwo{\left(\atsu X_s+\btsu X_u\right)\pln{n-1}{u}}{\mathcal{F}_{\le s}}=\\
&=&\atsu X_s \pln{n-1}{s}+\btsu \wwo{X_u\pln{n-1}{u}}{\mathcal{F}_{\le s}}=\\
&=&\atsu\big(a_{n-1}(s)\pln{n}{s}+b_{n-1}(s)\pln{n-1}{s}+c_{n-1}\pln{n-2}{s}\big)+\\
&+&\btsu\big(a_{n-1}(u)\wwo{\pln{n}{u}}{\mathcal{F}_{\le s}} + b_{n-1}(u)\pln{n-1}{s} +
c_{n-1}(u)\pln{n-2}{s}\big).
\end{eqnarray*}

On the other hand,

\begin{eqnarray*}
& &\wwo{X_t \pln{n-1}{u}}{\mathcal{F}_{\le s}}= \wwo{X_t\wwo{\pln{n-1}{u}}{\mathcal{F}_{\le
t}}}{\mathcal{F}_{\le s}}
=\wwo{X_t \pln{n-1}{t}}{\mathcal{F}_{\le s}}\\
&=&a_{n-1}(t)\wwo{\pln{n}{t}}{\mathcal{F}_{\le s}}+b_{n-1}(t)\pln{n-1}{s}+c_{n-1}(t)\pln{n-2}{s}.
\end{eqnarray*}

Thus  comparing the right hand sides of the above equations we obtain the following equation:

\[\btsu a_{n-1}(u)\wwo{\pln{n}{u}}{\mathcal{F}_{\le
s}}-a_{n-1}(t)\wwo{\pln{n}{t}}{\mathcal{F}_{\le s}}=\]
\[=-\pln{n}{s}\atsu a_{n-1}(s)+\pln{n-1}{s}\Big(b_{n-1}(t)-\atsu b_{n-1}(s)-\btsu b_{n-1}(u)\Big)\]
\[
+\pln{n-2}{s}\Big(c_{n-1}(t)-\atsu c_{n-1}(s)-\btsu c_{n-1}(u)\Big).
\]
A trivial verification, using \eqref{coeff}, shows that
\begin{eqnarray*}
b_{n-1}(t)-\atsu b_{n-1}(s)-\btsu b_{n-1}(u)&=&0,\\
c_{n-1}(t)-\atsu c_{n-1}(s)-\btsu c_{n-1}(u)&=&0.
\end{eqnarray*}
Hence denoting
\newcommand{\Pu}{\mathcal{X}}
\newcommand{\Pt}{\mathcal{Y}}

\[
\Pu:=\wwo{\pln{n}{u}}{\mathcal{F}_{\le s}},\ \Pt:=\wwo{\pln{n}{t}}{\mathcal{F}_{\le s}},
\]
we have
\begin{equation}\label{1eq}
\btsu a_{n-1}(u)\Pu-a_{n-1}(t)\Pt=-\atsu a_{n-1}(s)\pln{n}{s}.
\end{equation}

To obtain a second equation for $\Pu$ and $\Pt$ let us consider $\wwo{X_t^2
\pln{n-2}{u}}{\mathcal{F}_{\le s}}$.

On one hand
\begin{multline*}
\wwo{X_t^2 \pln{n-2}{u}}{\mathcal{F}_{\le s}}=\wwo{X_t^2\wwo{\pln{n-2}{u}}{\mathcal{F}_{\le
t}}}{\mathcal{F}_{\le s}}=\\
=\wwo{X_t^2 \pln{n-2}{t}}{\mathcal{F}_{\le s}}.
\end{multline*}
Setting $a_k(t)=b_k(t)=c_k(t)=0$ for $k<0$, a repeated application of \eqref{three step abstract} gives
\begin{multline}\label{2xthreeterm}
x^2 p_{n-2}\left(x;t\right)=p_n\left(x;t\right) a_{n-2}(t)a_{n-1}(t)
+p_{n-1}\left(x;t\right)a_{n-2}(t)\left[b_{n-2}(t)+b_{n-1}(t)\right]+\\
+p_{n-2}\left(x;t\right)\left[a_{n-2}(t)c_{n-1}(t)+b_{n-2}^2(t)+c_{n-2}(t)a_{n-3}(t)\right]+\\
+p_{n-3}\left(x;t\right)c_{n-2}(t)\left[b_{n-3}(t)+b_{n-2}(t)\right]+
p_{n-4}\left(x;t\right)c_{n-3}(t)c_{n-2}(t),
\end{multline}
we obtain
\begin{multline}\label{2eqL}
\wwo{X_t^2 \pln{n-2}{u}}{\mathcal{F}_{\le s}}=a_{n-2}(t)a_{n-1}(t)\ \Pt+\\
+\pln{n-1}{s}a_{n-2}(t)\left[b_{n-2}(t)+b_{n-1}(t)\right]+\\
+\pln{n-2}{s}\left[a_{n-2}(t)c_{n-1}(t)+b_{n-2}^2(t)+c_{n-2}(t)a_{n-3}(t)\right]+\\
+\pln{n-3}{s}c_{n-2}(t)\left[b_{n-3}(t)+b_{n-2}(t)\right]+\\
+\pln{n-4}{s}c_{n-3}(t)c_{n-2}(t).
\end{multline}
On the other hand, one can rewrite $\wwo{X_t^2 \pln{n-2}{u}}{\mathcal{F}_{\le s}}$ as
\[\wwo{X_t^2 \pln{n-2}{u}}{\mathcal{F}_{\le
s}}=\wwo{\wwo{X_t^2}{\mathcal{F}_{s,u}}\pln{n-2}{u}}{\mathcal{F}_{\le s}}\]
\[
=\Atsu X_s^2\pln{n-2}{s}+\Btsu X_s\wwo{X_u\pln{n-2}{u}}{\mathcal{F}_{\le
s}}\]
\[+
\Ctsu\wwo{X_u^2\pln{n-2}{u}}{\mathcal{F}_{\le s}}+ \Dtsu X_s\pln{n-2}{s}\]
\[
+\Etsu\wwo{X_u\pln{n-2}{u}}{\mathcal{F}_{\le s}}+ \Ftsu\pln{n-2}{s}.
\]
After some algebra one gets
\begin{multline}\label{Wojtek1}
\wwo{X_t^2 \pln{n-2}{u}}{\mathcal{F}_{\le s}}=\Pu\ \Ctsu a_{n-2}(u)a_{n-1}(u)
\\+\pln{n}{s}a_{n-1}(s)R_n(t,s,u)
+\pln{n-1}{s} S_n(t,s,u)
+\pln{n-2}{s}T_n(t,s,u)+\\
+\pln{n-3}{s}U_n(t,s,u)
+\pln{n-4}{s}V_n(t,s,u),
\end{multline}
where
\begin{eqnarray*}
R_n(t,s,u)&=&\Atsu a_{n-2}(s)+\Btsu a_{n-2}(u),\\
S_n(t,s,u)&=&\Atsu a_{n-2}(s)\left[b_{n-2}(s)+b_{n-1}(s)\right]
+\Btsu\left[a_{n-2}(u)b_{n-1}(s)+a_{n-2}(s)b_{n-2}(u)\right] \\
&&+\Ctsu a_{n-2}(u)\left[b_{n-2}(u)+b_{n-1}(u)\right] +\Dtsu b_{n-2}(s) +\Etsu
b_{n-2}(u)+\Ftsu\,, \\
T_n(t,s,u)&=&\Atsu\left[a_{n-2}(s)c_{n-1}(s)+b_{n-2}^2(s)+c_{n-2}(s)a_{n-3}(s)\right]\\
&&+\Btsu\left[a_{n-2}(u)c_{n-1}(s)+b_{n-2}(u)b_{n-2}(s)+c_{n-2}(u)a_{n-3}(s)\right]\\
&&+\Ctsu\left[a_{n-2}(u)c_{n-1}(u)+b_{n-2}^2(u)+c_{n-2}(u)a_{n-3}(u)\right] \\ &&
+\Dtsu b_{n-2}(s)+\Etsu b_{n-2}(u)+\Ftsu\,,\\
U_n(t,s,u)&=& \Atsu c_{n-2}(s)\left[b_{n-3}(s)+b_{n-2}(s)\right]+
\Btsu\left[b_{n-2}(u)c_{n-2}(s)+c_{n-2}(u)b_{n-3}(s)\right]\\
&&+\Ctsu c_{n-2}(u)\left[b_{n-3}(u)+b_{n-2}(u)\right]+ \Dtsu c_{n-2}(s) +\Etsu c_{n-2}(u),\\
V_n(t,s,u)&=& \Atsu c_{n-2}(s)c_{n-3}(s)+\Btsu c_{n-2}(u)c_{n-3}(s)+\Ctsu
c_{n-2}(u)c_{n-3}(u).
\end{eqnarray*}

Comparing the right hand side of \rf{Wojtek1} with the right hand side of \eqref{2eqL} we get
the second equation with unknowns $\Pu$ and $\Pt$:
\begin{multline}\label{2eq}
\Ctsu a_{n-2}(u)a_{n-1}(u)\ \Pu - a_{n-2}(t)a_{n-1}(t)\ \Pt=
-\pln{n}{s}R_n(t,s,u)\\+
\pln{n-1}{s} \Big(a_{n-2}(t)\left[b_{n-2}(t)+b_{n-1}(t)\right] -S_n(t,s,u)\Big)\\+
\pln{n-2}{s} \Big(a_{n-2}(t)c_{n-1}(t)+b_{n-2}^2(t)+c_{n-2}(t)a_{n-3}(t)-T_n(t,s,u)\Big)\\+
\pln{n-3}{s} \Big(c_{n-2}(t)\left[b_{n-3}(t)+b_{n-2}(t)\right] -U_n(t,s,u)\Big)\\+
\pln{n-4}{s} \Big(c_{n-2}(t) c_{n-3}(t) -V_n(t,s,u)\Big).
\end{multline}
A calculation based on \rf{coeff}  and   (\ref{A0}-\ref{E0}) gives
\begin{multline*}
a_{n-2}(t)\left[b_{n-2}(t)+b_{n-1}(t)\right] -S_n(t,s,u)\\
=-\frac{(t-s)(u-t)}{u(1+\sigma s)+\tau-qs} \Big[\theta\sigma\alpha_{n-1}+ \eta\beta_{n-1}
-\beta_{n-1}\gamma_{n-1} +\sigma\tau\alpha_{n-1}\left(\gamma_{n-2}+\gamma_{n-1}\right)-\\
-\sigma\alpha_{n-1}\delta_{n-2}+\sigma\beta_{n-1}\left(\delta_{n-2}+\delta_{n-1}\right)
+q\left(\beta_{n-1}\gamma_{n-2}+\sigma\alpha_{n-1}\delta_{n-1}\right)\Big],
\end{multline*}
\begin{multline*}
a_{n-2}(t)c_{n-1}(t)+b_{n-2}^2(t)+c_{n-2}(t)a_{n-3}(t)-T_n(t,s,u)\\
=-\frac{(t-s)(u-t)}{u(1+\sigma s)+\tau-qs}\Big[ 1 + \theta\gamma_{n-2}+ \eta\delta_{n-2}
-\gamma_{n-2}\delta_{n-2} -\beta_{n-1}\epsilon_{n-1}+\\
+\tau\left({\gamma}_{n-2}^2
+\sigma\alpha_{n-2}\epsilon_{n-2}+\sigma\alpha_{n-1}{\eps}_{n-1}\right)-\sigma\alpha_{n-2}\phi_{n-2}+\\
+q\left(\gamma_{n-2}\delta_{n-2}+\beta_{n-2}\epsilon_{n-2}+ \sigma\alpha_{n-1}\phi_{n-1}\right)  +
\sigma\left(\delta_{n-2}^2 +\beta_{n-2}\phi_{n-2} + \beta_{n-1}\phi_{n-1}\right)\Big],
\end{multline*}
\begin{multline*}
c_{n-2}(t)\left[b_{n-3}(t)+b_{n-2}(t)\right] -U_n(t,s,u)=
-\frac{(t-s)(u-t)}{u(1+\sigma s)+\tau-qs}
\Big[\theta\eps_{n-2}\\+\tau\eps_{n-2}\left(\gamma_{n-3}+\gamma_{n-2}\right)-\delta_{n-2}\eps_{n-2}+
\eta\phi_{n-2}-\phi_{n-2}\gamma_{n-3}\\
+\sigma\left(\delta_{n-3}+\delta_{n-2}\right)\phi_{n-2}
+q\left(\delta_{n-3}\eps_{n-2}+\gamma_{n-2}\phi_{n-2}\right) \Big],
\end{multline*}
\begin{multline*}
c_{n-2}(t) c_{n-3}(t) -V_n(t,s,u)\\\
=-\frac{(t-s)(u-t)}{u(1+\sigma s)+\tau-qs} \big[\eps_{n-3}\left(\tau\eps_{n-2}-\phi_{n-2}\right)+
\phi_{n-3}\left(q\eps_{n-2}+\sigma\phi_{n-2}\right)\big],
\end{multline*}
Therefore, by \eqref{eq2}-\eqref{eq5} the coefficients at $\pln{n-1}{s},\ldots,\pln{n-4}{s}$ on the right hand side
of \eqref{2eq} vanish, and \eqref{2eq} is equivalent to
\begin{multline}\label{2eqmod}
\Ctsu a_{n-2}(u)a_{n-1}(u)\ \Pu - a_{n-2}(t)a_{n-1}(t)\ \Pt \\
=
- a_{n-1}(s)\pln{n}{s}\left(\Atsu a_{n-2}(s)+\Btsu a_{n-2}(u)\right).
\end{multline}
Subtracting from \rf{2eqmod} equation \rf{1eq} multiplied by $a_{n-2}(t)$, we get
\begin{multline}\label{2eq ver3}
a_{n-1}(u)\left(\Ctsu a_{n-2}(u)-\btsu a_{n-2}(t)\right)\Pu \\
=
- a_{n-1}(s)\pln{n}{s}\left(\Atsu a_{n-2}(s)+\Btsu a_{n-2}(u)-\atsu a_{n-2}(t)\right).
\end{multline}

Now we notice that
\begin{multline*}
\Ctsu a_{n-2}(u)-\btsu a_{n-2}(t)=\\
=- \frac{(t-s)(u-t)}{(u-s)\left[u(1+\sigma s)+\tau-q s\right]} \left[(q
s-\tau)\sigma\alpha_{n-1}+(1+s\sigma)\beta_{n-1}\right]
\end{multline*}
see \rf{non-linear rec} and \rf{F0}.
Similarly,
\begin{multline*}\Atsu a_{n-2}(s)+\Btsu a_{n-2}(u)-\atsu a_{n-2}(t)= \\
=
\frac{(t-s)(u-t)}{(u-s)\left[u(1+\sigma s)+\tau-q s\right]} \left[(q
u-\tau)\sigma\alpha_{n-1}+(1+u\sigma)\beta_{n-1}\right] 
\end{multline*}
Since \rf{eq1} implies that
\begin{multline}\label{wstawka}
a_{n-1}(u)\left((q
s-\tau)\sigma\alpha_{n-1}+(1+s\sigma)\beta_{n-1}\right)\\ =a_{n-1}(s)\left((q
u-\tau)\sigma\alpha_{n-1}+(1+u\sigma)\beta_{n-1}\right),
\end{multline}
and $a_{n-1}(u)\ne 0$ by assumption, therefore, \rf{2eq ver3} becomes
\begin{equation}\label{2eq ver4}
\left((q
s-\tau)\sigma\alpha_{n-1}+(1+s\sigma)\beta_{n-1}\right)\Pu =
\left((q
s-\tau)\sigma\alpha_{n-1}+(1+s\sigma)\beta_{n-1}\right)\pln{n}{s}.
\end{equation}
It remains to verify that
\begin{equation}\label{wstawka2}
\left((q
s-\tau)\sigma\alpha_{n-1}+(1+s\sigma)\beta_{n-1}\right)\ne 0.
\end{equation}
Since  $\beta_{n-1}=a_{n-1}(0)\ne 0$ by assumption, this is trivially true if $\sigma=0$.
Suppose $\sigma>0$ and \rf{wstawka2} isn't true. Since $a_{n-1}(u)\ne 0$
from \rf{wstawka} we see that
then the left hand side of \rf{wstawka2} must vanish also when $s$ is replaced by $u$.
Thus
\begin{equation}\label{tmptmp1}
\sigma(\alpha_{n-1} q+\beta_{n-1})=0, \; \beta_{n-1}-\sigma\tau\alpha_{n-1}=0.
\end{equation}
Since $\beta_{n-1}\ne 0$,
the determinant of the system of equations \eqref{tmptmp1}  is zero,
 $q+\sigma\tau=0$, and the second equation gives $\sigma\tau\ne 0$.
Since $\alpha_1=0$, from \rf{eq1} with $q=-\sigma\tau\ne 0$ we get $\alpha_j=\beta_j$ for all $j\geq 1$. As $\sigma\tau\ne 1$,
this  contradicts \eqref{tmptmp1}.
Thus \rf{wstawka2} holds.

By \rf{wstawka2}, from \rf{2eq ver4} we get
\[
\Pu=\pln{n}{s},
\]
which means that $\left\{p_n(x;t)\right\}_{t>0}$ is a martingale for every $n\geq 0$.

Now we shall show that polynomials $\left\{p_n(x;t)\right\}_n$ are orthogonal with respect to the
distribution of $X_t$. Since the polynomials  satisfy the three term recurrence
\eqref{three step abstract}, it suffices to show that for all $k=1,2,\ldots$
\begin{equation}\label{ompind2}
\wo{(\pln{k}{t})}=0.
\end{equation}
 Clearly, by the martingale property, $\mu_k=\wo{(\pln{k}{t})}$ does not depend on $t$.

We proceed again by induction. From  \rf{non-linear rec} we see that $a_1(t)=1+\sigma t$ so
$(1+\sigma t) p_2(x;t)=x^2-b_1(t)p_1(x;t)-t$. Since $p_1(x;t)=x$,  assumption \rf{cov} implies that \eqref{ompind2} holds for $k=1,2$.

Assume then that for some $n\geq 3$ we have $\mu_k=0$ for all $1\leq k\leq n-1$. By
\eqref{three step abstract}
\[
\wo{(X_t\pln{n}{t})}=a_n(t)\wo{(\pln{n+1}{t})}=a_n(t)\mu_{n+1}.
\]
By \eqref{LR+}, the martingale property of $\left(p_n(x;t)\right)_n$, \eqref{three step abstract}
and the induction assumption, the left hand side of the above equation transforms into
\begin{multline*}
\wo{(X_t\pln{n}{t})}=\wo\Big(\frac{t}{s}{\wwo{X_s}{\mathcal{F}_{\ge
t}}\pln{n}{t}}\Big)=\frac{t}{s}\wo\Big(X_s\wwo{\pln{n}{t}}{\mathcal{F}_{\le s}}\Big)=\\
=\frac{t}{s}\wo\Big(X_s\pln{n}{s}\Big)=\frac{t}{s}a_n(s)\wo{(\pln{n+1}{s})}= \frac{t}{s}a_n(s)\mu_{n+1},
\end{multline*}
so
\[
\left(\frac{a_n(s)}{s}-\frac{a_n(t)}{t}\right)\mu_{n+1}=0.
\]
Since  $$\frac{a_n(s)}{s}-\frac{a_n(t)}{t}=a_n(0)\frac{t-s}{ts}\ne 0,$$
 we get
 $\mu_{n+1}=\wo{(\pln{n+1}{t})}=0$, and the orthogonality of
$\left\{p_n(x;t)\right\}_n$ is proved.


%


\subsection{Proof of Theorem \ref{T4.1}}
\begin{lemma}\label{BP-tail}
Suppose $X,Y$ are square-integrable standardized,  and
there are constants $\eps,A,B\geq 0, 0<\rho<1$ such that
\begin{eqnarray}
\wwo{(X-\rho Y)^2}{Y}&\leq& A+ B|Y| +\frac{(1-\rho^2)\eps^2}{1+\eps^2}Y^2, \label{BP1}\\
\wwo{(Y-\rho X)^2}{X}&\leq& A+ B |X| +\frac{(1-\rho^2)\eps^2}{1+\eps^2}X^2. \label{BP2}
\end{eqnarray}
If $2\eps<\rho<1$, then  there are constants $C_1,C_2<\infty $ such that for $x>0$ we have
\begin{eqnarray} \label{BP-tail bound}
&&\Pr(|X|>2x/\rho)+\Pr(|Y|>2x/\rho)  \\ \nonumber
&&\leq \left(C_1/x^2+C_2/x+\frac{4\eps^2}{\rho^2+\eps^2}\frac{1+\rho}{1-\rho}\right)\left(\Pr(|X|>x)+\Pr(|Y|>x)\right).
\end{eqnarray}
\end{lemma}

\begin{proof} Enlarging $\eps$ if necessary, without loss of generality we may assume that $\eps>0$.

Let $N(x)=P(|X|\geq x)+P(|Y|\geq x)$, $K=2/\rho$. 
Throughout the proof, $C_1,C_2$ denote positive constants which might
 differ at each occurrence.

The event $\{|X|\geq Kx\}$, where $x>0$ is fixed,
can be decomposed into the sum of two disjoint events
$\{|X|\geq Kx\}\cap\{|Y|\geq x\}$
 and $\{|X|\geq Kx\}\cap \{|Y|<x\}$. Therefore,
denoting
\begin{eqnarray*}
P_1(x)&=&\Pr(|X|\geq x, |Y|\geq x),\\
P_2(x)&=&\Pr(|X|\geq Kx, |Y|<x),\\
 P_3(x)&=&\Pr(|Y|\geq Kx, |X|<x),
\end{eqnarray*} we get
$
 \Pr(|X|\geq Kx)\leq P_1(x)+P_2(x),
 $
 and hence by symmetry of the assumptions
 \begin{equation}
\label{(4.3.8)}
N(Kx)\leq 2P_1(x)+P_2(x)+P_3(x).
\end{equation}

To estimate $P_1(x)$, we  observe the following. 
\begin{claim}\label{Cl 3}
  If $\{|X|\geq x$, $|Y|\geq x\}$, and $0<a<1$, then
either
$$(X-\rho Y)^2\geq a(1-\rho)^2 Y^2+(1-\rho)^2(1-a)x^2,$$
 or
 $$(Y-\rho X)^2\geq a (1-\rho)^2X^2 + (1-\rho)^2(1 -a)x^2.$$
 \end{claim}
Indeed, suppose that both inequalities fail, i.e. on the set $\{|X|\geq x$, $|Y|\geq x\}$ we have
$(X-\rho Y)^2< a(1-\rho)^2 X^2+(1-\rho)^2(1-a)x^2$ and
$(Y-\rho X)^2< a(1-\rho)^2Y^2 + (1-\rho)^2(1 -a)x^2$.
Adding the inequalities we obtain
$$(1+\rho^2)(X^2+Y^2)-4\rho XY< a(1-\rho)^2 (X^2+Y^2)+2 (1-\rho)^2(1 -a )x^2.$$
Since $4\rho XY\leq 2 \rho(X^2+Y^2)$, this gives
\begin{equation*}\label{12345}
(1-\rho)^2(1-a)(X^2+Y^2)< 2 (1-\rho)^2(1 -a )x^2.
\end{equation*}
However, since $|X|\geq x,|Y|\geq x$, we have
$$(1-\rho)^2(1-a)(X^2+Y^2) \geq 2 x^2(1-\rho)^2(1-a ),
$$
a contradiction.

Claim \ref{Cl 3}  with $a=1/2$  implies
$$P_1(x)\leq \Pr(|X-\rho Y|\geq (1-\rho)\sqrt{Y^2+x^2}/\sqrt{2}, |Y|\geq x)$$
$$+\Pr(|Y-\rho X|\geq (1-\rho)\sqrt{ X^2+x^2}/\sqrt{2}, |X|\geq x).$$
From conditional Chebyshev's inequality and \rf{BP1} we get
$$\Pr(|X-\rho Y|\geq (1-\rho)\sqrt{Y^2+x^2}/\sqrt{2}, |Y|\geq x)$$
$$
\leq 2\int_{|Y|\geq x}\frac{A+ B|Y| +(1-\rho^2){\eps^2}Y^2/{(1+\eps^2)}}{(1-\rho)^2( Y^2+x^2)}dP
$$
$$
\leq
\int_{|Y|\geq x}\frac{2A}{(1-\rho)^2x^2}dP+
\int_{|Y|\geq x}\frac{2B|Y|}{(1-\rho)^2 |Y| x}dP+
\int_{|Y|\geq x}\frac{(1-\rho^2){2\eps^2}Y^2/{(1+\eps^2)}}{(1-\rho)^2Y^2}dP
$$
$$
\leq \frac{C_1}{x^2} \Pr(|Y|\geq x)+\frac{C_2}{x} \Pr(|Y|\geq x)+
\frac{2\eps^2}{1+\eps^2}\frac{1+\rho}{1-\rho} \Pr(|Y|\geq x),
$$
where $C_1=2A/(1-\rho)^2$, $C_2=2B/(1-\rho)^2$.
Since the assumptions are symmetric in $X,Y$, this shows that there are constants $C_1,C_2<\infty$ such that

\begin{equation}\label{21345}
P_1(x)
\leq \frac{C_1 N(x)}{x^2}+\frac{C_2 N(x)}{x}+
\frac{2\eps^2}{1+\eps^2}\frac{1+\rho}{1-\rho} N(x).
\end{equation}
To estimate $P_2(x)$, we use the trivial estimate $|Y-\rho X|\geq \rho |X|-|Y|$, which shows that the event
$\{|X|\geq Kx, |Y|<x\}$ implies
\[ |Y-\rho X|\geq \eps |X|+\left((\rho-\eps)K-1\right)x=\eps |X|+\frac{\rho-2\eps}{\rho}x
\]
Therefore
\begin{align*}
P_2(x) \leq & \Pr(|Y-\rho X|\geq \eps |X|+\frac{\rho-2\eps}{\rho}x, |X|\geq Kx)& \\ \nonumber
&=\int_{|X|\geq Kx}\Pr(|Y-\rho X|\geq \eps |X|+\frac{\rho-2\eps}{\rho}x|X)\, dP \\ \nonumber
&\displaystyle \leq\int_{|X|\geq Kx} \frac{A+ B |X| +(1-\rho^2){\eps^2}X^2/{(1+\eps^2)}}{(\eps |X|+{(\rho-2\eps)}x/{\rho})^2} dP &
\\ \nonumber
& \displaystyle\leq
\int_{|X|\geq Kx} \frac{A}{((\rho-2\eps)x/{\rho})^2} dP
+\int_{|X|\geq Kx} \frac{B |X|}{\eps |X|(\rho-2\eps)x/\rho} dP&\\
&
+\int_{|X|\geq Kx} \frac{(1-\rho^2)\eps^2X^2}{(1+\eps^2)\eps^2 X^2} dP.
&
 \end{align*}
 Since $N(Kx)\leq N(x)$, this shows that
 \begin{equation}\label{(4.3.9)}
P_2(x) \leq \frac{C_1 N(x)}{x^2}+\frac{C_2N(x)}{x}+\frac{1-\rho^2}{1+\eps^2}\Pr(|X|\geq Kx).
 \end{equation}
where $C_1=\rho^2A/(\rho-2\eps)^2$, $C_2=\rho B/(\eps(\rho-2\eps))$.
By symmetry of the assumptions, we also have
\begin{equation}\label{(4.3.9')}
P_3(x) \leq \frac{C_1 N(x)}{x^2}+\frac{C_2N(x)}{x}+\frac{1-\rho^2}{1+\eps^2}\Pr(|Y|\geq Kx).
 \end{equation}

Combining (\ref{(4.3.8)}), \rf{21345}, (\ref{(4.3.9)}), and (\ref{(4.3.9')}) we obtain that there are constants
$C_3,C_4>0$ such that
$$
N(Kx)\leq  \frac{C_3 N(x)}{x^2}+\frac{C_4N(x)}{x}+\frac{1-\rho^2}{1+\eps^2} N(Kx)
+ \frac{4 \eps^2}{1+\eps^2}\frac{1+\rho}{1-\rho} N(x),
$$
which implies \rf{BP-tail bound}.
\end{proof}

\begin{lemma} \label{BP integr}Suppose $X,Y$ are square-integrable standardized,
$\wo(|X|^p)+\wo(|Y|^p)<\infty$ for some $p> 1$ and
the assumptions of Lemma  \ref{BP-tail} are satisfied with constants $\eps> 0, 0<\rho<1$ such that
 \begin{equation}\label{BP bound}
2^{p+3} \eps^2\frac{1+\rho}{1-\rho}\leq \rho^{p+3}.
\end{equation}
Then $\wo(|X|^{p+1})+\wo(|Y|^{p+1})<\infty$.
\end{lemma}
\begin{proof} Clearly, \rf{BP bound} implies $2\eps<\rho$. Indeed,
$$2^2\eps^2<2^{p+3}\eps^2\frac{1+\rho}{1-\rho}\leq \rho^{p+3}<\rho^2.$$ We use Lemma \ref{BP-tail} and we use
the notation $N(x)=P(|X|\geq x)+P(|Y|\geq x)$, $K=2/\rho$ introduced in its proof.
Fix $M>0$. Then, noticing that $K>1$ from Lemma \ref{BP-tail} we get
$$
\int_0^M (p+1) x^p N(x)dx= K^{p+1}\int_0^{M/K} (p+1) x^p N(Kx)dx
$$
$$
\leq \frac{p+1}{p-1}C_1 \int_0^\infty (p-1)x^{p-2} N(x) dx+ \frac{p+1}{p}C_2 \int_0^\infty px^{p-1} N(x) dx$$
$$+
\frac{4 \eps^2}{\rho^2+\eps^2}\frac{1+\rho}{1-\rho} K^{p+1}\int_0^M (p+1)x^{p} N(x) dx
$$
$$<\frac{p+1}{p-1}C_1(\wo(|X|^{p-1})+\wo(|Y|^{p-1}))+\frac{p+1}{p}C_2(\wo(|X|^{p})+\wo(|Y|^{p}))$$
$$+\eps^2 2^{p+3}/\rho^{p+3}\frac{1+\rho}{1-\rho}\int_0^M (p+1)x^{p} N(x) dx.
$$
Therefore, if \rf{BP bound} holds true, then
$\sup_{M>0}\int_0^M (p+1) x^p N(x)dx<\infty$ which implies
$\wo(|X|^{p+1})+\wo(|Y|^{p+1})=\int_0^\infty (p+1) x^p N(x)dx<\infty$.
\end{proof}
 We will also need the following simple observation.
\begin{lemma} \label{claim 4}
Fix $r\geq 1$. Under the assumptions of Theorem \ref{T4.1}, if $\wo(|X_t|^r)<\infty $ for some $t>0$, then
$\wo(|X_t|^{r})<\infty$ for all $t>0$.
\end{lemma} \begin{proof}
If $\wo(|X_t|^{r})<\infty$ for some $t>0$ then
 $\wo(|X_s|^{r})<\infty$ for all $s<t$ by \rf{LR-}. Similarly,
$\wo(|X_u|^{r})<\infty$ for all $u>t$
by \rf{LR+}.
\end{proof}

\begin{proof}[Conclusion of proof of Theorem \ref{T4.1}]

Let $r> 2$  be  fixed. Since our proof relies only on inequalities, if $\sigma\tau=0$
we can increase slightly the values of parameters $\sigma, \tau$ on the right hand sides of
\rf{RV-} and \rf{CV-} to ensure $\sigma',\tau'>0$ and $\sigma'\tau'\leq 2^{-4r-10}$,
replacing both equalities by the appropriate inequalities;
 thus without loss of generality we may assume that  $\sigma\tau>0$.

Let $p=r-1> 1$.
Define
$\eps^2=\sqrt{\sigma\tau}/2$, $t_0=\tau/\eps^2, s_0=\eps^2/\sigma$.
Then the correlation coefficient between $X_{s_0}$ and $X_{t_0}$ is $\rho=\sqrt{s_0/t_0}=\eps^2/\sqrt{\sigma\tau}=1/2$.
From the inequality form of \rf{RV-} and \rf{CV-} 
we deduce that
random variables $X=X_{s_0}/\sqrt{s_0}$ and $Y=X_{t_0}/\sqrt{t_0}$ satisfy \rf{BP1} and \rf{BP2} with
$$
A=\frac{1-\rho^2}{1+\eps^2}, \, B=\frac{1-\rho^2}{1+\eps^2}\max \{|\eta|\eps/\sqrt{\sigma},|\theta|\sqrt{\tau}/\eps\}.
$$
Condition \rf{BP bound} is satisfied since with $r=p+1$ we have
$(\sigma \tau)^{1/2}\leq 1/2^{2p+7}$. Thus
$$
2^{p+3} \eps^2\frac{1+\rho}{1-\rho}=
2^{p+3} 3\eps^2< 2^{p+4}(\sigma \tau)^{1/2}\leq 1/2^{p+3}= \rho^{p+3}.
$$
Therefore, by Lemma \ref{BP integr}, $\wo(|X_{t_0}|^r)=\wo(|X_{t_0}|^{p+1})<\infty$ and therefore
by Lemma \ref{claim 4},
$\wo(|X_t|^{r})<\infty$ for all $t\geq 0$.
\end{proof}

\section{Examples of explicit recurrences}\label{SECT: Appl}
In this section we will say that $(X_t)$ is a quadratic harness with parameters $q,\eta,\theta,\sigma,\tau$
if \rf{LH}, \rf{cov},  and \rf{q-Var} hold.
Our goal is to
derive explicit versions of \rf{three step abstract}  at the expense of additional assumptions on the parameters.


\begin{theorem}
\label{Thm C-com+}
 Suppose $(X_t)$ is a quadratic harness
with parameters
$0\leq \sigma\tau<1$, $-1< q\leq 1-2\sqrt{\sigma\tau}$.
Moreover, assume that for each $t>0$ random variable $X_t$  has 
moments of all orders and infinite support.
Then $(X_t)$ has orthogonal martingale polynomials which are given by recurrence \rf{three step abstract}
with  $p_0(x,t)=1$, $p_1(x,t)=x$, where the coefficients of the recurrence are linear functions of $t$,
%
\begin{equation}\label{coeff+}
a_n(t)=\sigma\alpha_{n+1} t+\beta_{n+1}\;,\;\;\;b_n(t)=\gamma_nt+\delta_n\;,\;\;
\;c_n(t)=\left(\beta_nt+\tau\alpha_n\right)\omega_n\,,
\end{equation}
which are determined as follows.

\begin{enumerate}

\item The initial conditions are
$$\alpha_1=0,\; \beta_1=1,\;\gamma_0=\delta_0=0,\;\omega_1=1.$$

\item Sequences
$(\alpha_n),\;(\beta_n)$ satisfy 
\begin{equation}\label{non-linear rec}
\left[
\begin{array}{c}
  \alpha_{n+1} \\
  \beta_{n+1} \\
\end{array}
\right] =
\left[
\begin{array}{cc}
  q & 1 \\
  -\sigma \tau & 1 \\
\end{array}
\right]\times
\left[
\begin{array}{c}
  \alpha_n \\
  \beta_n \\
\end{array}
\right], \quad n\geq 1.
\end{equation}
 Moreover,
denoting $$\la_{n,k}=\beta_{n }\beta_{n+k }-\sigma\tau \alpha_{n }\alpha_{n+k }$$
we have $\la_{n,k}>0$ for all $n\geq 1$, $k\geq 0$.
\item Setting $\la_{0,2}=\alpha_0=\beta_0=0$, sequences
 $(\gamma_n),\;(\delta_n)$ satisfy the following system of linear recurrences
 \begin{eqnarray}\label{gamma-delta rec} \nonumber
\gamma_{n+1}&=& \frac{q+\sigma\tau}{\la_{n+2,0}}\left(\la_{n,2}\gamma_n+
(\alpha_{n+2}\beta_n- \beta_{n+2}\alpha_n)\sigma\delta_n\right) \\
       &&  +
         \frac{\sigma\alpha_{n+2}}{\la_{n+2,0}}(\eta\tau\alpha_{n+1}+\theta\beta_{n+1})+
         \frac{\beta_{n+2}}{\la_{n+2,0}}(\theta\sigma\alpha_{n+1}+\eta\beta_{n+1}),
\\ \nonumber
\delta_{n+1}&=&
\frac{q+\sigma\tau}{\la_{n+2,0}}\left(\la_{n,2}\delta_n+
(\alpha_{n+2}\beta_n- \beta_{n+2}\alpha_n)\tau\gamma_n\right)
\\ \nonumber&&  +
         \frac{\beta_{n+2}}{\la_{n+2,0}}(\eta\tau\alpha_{n+1}+\theta\beta_{n+1})+
         \frac{\tau\alpha_{n+2}}{\la_{n+2,0}}(\theta\sigma\alpha_{n+1}+\eta\beta_{n+1}), \; n\geq 0.
\end{eqnarray}

\item Sequence $(\omega_n)$ satisfies the linear recurrence
\begin{eqnarray}\label{w rec}
\omega_{n+1}&=&\frac{(q+\sigma\tau)\la_{n-1,1}}{\la_{n+1,1}}\omega_n
\\ \nonumber
&&+
\frac{1+\gamma_n(\tau\gamma_n-\delta_n+\theta)+
\delta_n(q\gamma_n+\sigma\delta_n+\eta)}{\la_{n+1,1}},\; n\geq 2
\end{eqnarray}


with the initial term
$$ \omega_2=(1+q) \frac{( 1 - \sigma  \tau )^2 +(\eta+\theta\sigma)(\theta+\eta\tau)}
{(1-\sigma\tau)^2(1-\sigma\tau(2+q))}.$$
\end{enumerate}
\end{theorem}

\begin{remark} \label{R: Explicit solution}
If 
$-1\leq q< 1-2\sqrt{\sigma\tau}$, then
the explicit solution of  recurrence \rf{non-linear rec} is
\begin{equation}\label{alpha-beta}
\alpha_n=\frac{\la_+^{n-1}-\la_-^{n-1}}{\la_+-\la_-},\;
\beta_n=\frac{\la_+^{n-1}(1-\la_-)+\la_-^{n-1}(\la_+-1)}{\la_+-\la_-}
,
\end{equation}
where $$\la_\pm=\frac12\left(1+q\pm \sqrt{(1+q)^2-4(q+\sigma\tau)}\right).$$
\end{remark}

The proof of Theorem \ref{Thm C-com+} relies on Theorem \ref{C mo poly} and the
following lemma.

\begin{lemma}\label{L G3} 
If $\sigma\tau<1$ and  $-1\leq q\leq 1-2\sqrt{\sigma\tau}$, then the solutions of \rf{non-linear rec} satisfy
$$\beta_n>\sqrt{\sigma\tau}\alpha_n\geq 0.$$
In particular,
$$
\beta_n\beta_{n+1} >\sigma\tau \alpha_n \alpha_{n+1},\; \beta_n^2 >\sigma\tau \alpha_n^2\;, \mbox{ for all } n\geq 1.
$$
\end{lemma}

\begin{proof}From the explicit solution in Remark \ref{R: Explicit solution},
  we see that $|\la_-|\leq |\la_+|$ and hence
$\alpha_{n}\geq 0$ for $n>1$.
We now prove  by induction that $\beta_n>\sqrt{\sigma\tau}\alpha_n$ for all $n$.
This is trivially true for $n=1$, see \rf{ini}. Suppose the inequality is satisfied for some $n\geq 1$.
Then
\rf{non-linear rec} implies that
$\beta_{n+1}-\sqrt{\sigma\tau}\alpha_{n+1}=(1-\sqrt{\sigma\tau})\beta_n-\sqrt{\sigma\tau}(q+\sqrt{\sigma\tau})\alpha_n \geq
(1-\sqrt{\sigma\tau})(\beta_{n}-\sqrt{\sigma\tau}\alpha_{n})>0$.
\end{proof}

\begin{proof}[Proof of Theorem \ref{Thm C-com+}]

By Lemma \ref{L G3}, $a_n(t)=\sigma\alpha_{n+1} t+\beta_{n+1}\geq \beta_{n+1}>0$ for all $t\geq 0$, $n\geq 0$. Moreover,
recurrences \rf{gamma-delta rec}, \rf{w rec} are well defined and have unique solution.
Therefore, to end the proof we need only to verify that with
 \begin{equation}\label{www}
\eps_n=\omega_n\beta_n, \; \varphi_n= {\tau} \omega_n\alpha_n,
\end{equation}
the equations (\ref{ini}-\ref{eq5}) are satisfied. Trivially, \rf{ini} holds true.
It is easy to verify that the solution of  \rf{non-linear rec} satisfies \rf{eq1}. Using
\rf{www} from \rf{non-linear rec} we also get \rf{eq5}.

It remains to show that equations \rf{eq2}, \rf{eq3}, and \rf{eq4} hold. We will rewrite these equations using
 \rf{non-linear rec} and its equivalent form
\begin{equation}\label{Inverse}
q\beta_{n+1}+\sigma\tau\alpha_{n+1}= (q+\sigma\tau) \beta_n,\; \beta_{n+1} - \alpha_{n+1}=-(q+\sigma\tau)\alpha_n, \;n\geq1.
\end{equation}

Equation \rf{eq2}
is equivalent to
$$\gamma_{n+1}(\beta_{n+1}-\sigma\tau\alpha_{n+1})-\delta_{n+1}\sigma(\beta_{n+1}+\alpha_{n+1} q)
$$
$$=(\sigma\tau\alpha_{n+1}+\beta_{n+1}q)\gamma_n
+\delta_n\sigma(\beta_{n+1}-\alpha_{n+1})
+\sigma\alpha_{n+1}\theta+\beta_{n+1}\eta.$$
Using \rf{non-linear rec} and  \rf{Inverse}, we can rewrite it as
\begin{multline}\label{eq2'}
\gamma_{n+1}\beta_{n+2}-\sigma\delta_{n+1}\alpha_{n+2}
\\
=(q+\sigma\tau)\beta_{n}\gamma_n
-(q+\sigma\tau)\sigma\alpha_n\delta_n
+\sigma\alpha_{n+1}\theta+\beta_{n+1}\eta.
\end{multline}

Similarly, \rf{eq4} for our sequences rewrites as
\begin{multline}\label{eq2''}
\omega_n\left(\delta_{n+1}(\beta_{n+1}-\sigma\tau \alpha_{n+1})-\tau\gamma_{n+1}(\beta_{n+1}+q\alpha_{n+1})\right)
\\
=
\omega_n\left((q \beta_{n+1}+\sigma\tau\alpha_{n+1})\delta_n+{\tau}(\beta_{n+1}-\alpha_{n+1})\gamma_n+\theta\beta_{n+1}+
{\eta\tau}\alpha_{n+1}\right).
\end{multline}
Both of these equations are satisfied; in fact, \rf{gamma-delta rec} was obtained by solving the system of equations
(\ref{eq2'}-\ref{eq2''}) when $\omega_n>0$.

To verify that \rf{eq3} holds, we substitute \rf{www};
 using \rf{non-linear rec} and Lemma \ref{L G3} we get
\begin{multline*}
\omega_{n+1}=
\frac{\beta_n(q\beta_n+\sigma\tau\alpha_n)+
\sigma\tau\alpha_n(\beta_n - \alpha_n)}{\beta_{n+1}\beta_{n+2}-\sigma\tau\alpha_{n+1}\alpha_{n+2}} \omega_n
\\ +
\frac{1+\gamma_n(\tau\gamma_n-\delta_n+\theta)+
\delta_n(q\gamma_n+\sigma\delta_n+\eta)}{\beta_{n+1}\beta_{n+2}-\sigma\tau\alpha_{n+1}\alpha_{n+2}}.
\end{multline*}
Since \rf{Inverse} holds, thus \rf{eq3} is satisfied, too.
\end{proof}

\begin{remark}
From
\rf{gamma-delta rec} we get
\begin{equation}\label{delta1 gamma1}
\gamma_1=\frac{\eta+\theta\sigma}{1-\sigma\tau},\,\delta_1=\frac{\eta\tau+\theta}{1-\sigma\tau},
\end{equation}
which implies that
the initial recurrences
are \begin{align}
x p_0(x;t)&=p_1(x;t), \label{generic 0} \\
xp_1(x;t)&= (1+\sigma t)p_2(x;t)+  \frac{(\eta+\theta\sigma)t+(\eta\tau+\theta)}{1-\sigma\tau}p_1(x;t) + t p_0(x;t). \label{generic 1}
\end{align}
Thus the first three orthogonal martingale
 polynomials are
\begin{multline}\label{low order ortho}
p_0(x;t)=1,\;  p_1(x;t)= x, \\
p_2(x;t)=\frac{1}{1+\sigma t}x^2 -\frac{(\eta +\theta \sigma)t+\eta\tau+\theta}{(1-\sigma\tau)(1+\sigma t)}x-\frac{t}{1+\sigma t}.
\end{multline}
\end{remark}


\subsection{Free quadratic harnesses} Free harnesses have parameter
\begin{equation}\label{super-free}
q=-\sigma\tau.
\end{equation}
The adjective "free"  is motivated by the fact
that when $\sigma\tau=0$  this choice is  related to free convolutions that arise in free probability,
see \cite[Section 4.3]{Bryc-Wesolowski-03} and \cite[Section 4]{Bryc-Wesolowski-04}.
In general,
from Theorem \ref{Thm C-com+} it is clear that  this choice of $q$ significantly simplifies the recurrences.
It is not obvious whether this case is further related to free convolutions, see however \cite{Krystek-Wojakowski04}.


Assuming $\sigma\tau<1$,  \rf{super-free} implies $-1<q<1-2\sqrt{\sigma\tau}$, so it follows from Theorem \ref{Thm C-com+}
that the orthogonal martingale polynomials exist.
It is easy to check that
 the solution of \rf{non-linear rec} is
\begin{equation}\label{alpha-beta0}
\alpha_n=(1-\sigma\tau)^{n-2},\; \beta_n=(1-\sigma\tau)^{n-2},\;n\geq 2.
\end{equation}

For our choice of $q$ and $n\geq 0$ we have

$$
\la_{n+2,0}=\beta_{n+2}^2-\sigma\tau\alpha_{n+2}^2=(1-\sigma\tau)^{2n+1}.
$$

Using this identity, from
\rf{gamma-delta rec} we get  for $n\geq 2$
$$
\gamma_n=\frac{\eta  + 2\,\theta \,\sigma  + \eta \,\sigma \,\tau }
  {{\left( 1 - \sigma \,\tau  \right) }^2},
$$
$$
\delta_n=\frac{\theta  + 2\,\eta \,\tau  + \theta \,\sigma \,\tau }
  {{\left( 1 - \sigma \,\tau  \right) }^2}.
$$

Similarly,
$$
\la_{n+1,1}=\beta_{n+1}\beta_{n+2}-\sigma\tau\alpha_{n+1}\alpha_{n+2}=(1-\sigma\tau)^{2n}, \; n\geq 0
$$

Substituting  this and $\delta_n,\gamma_n$ into \rf{w rec}, we get
$$
\omega_n=\frac{ ( 1 - \sigma  \tau)^2 +
    (\eta+\theta\sigma)(\theta+\eta\tau) }
    {{\left( 1 - \sigma  \tau  \right) }^{2n}},\; n\geq 3.
$$ We also compute
$$
\omega_2=\frac{( 1 - \sigma  \tau )^2 +(\eta+\theta\sigma)(\theta+\eta\tau) }
    {{\left( 1 - \sigma  \tau  \right) }^{3}}.
$$

Therefore recurrence \rf{three step abstract} with initial values \rf{ini}
 gives
$$x p_2(x)=(1-\sigma\tau)(1+\sigma t)p_3(x)+
\frac{(\eta  + 2\,\theta \sigma  + \eta \,\sigma \,\tau )t+\theta  + 2\eta \tau  + \theta \sigma \tau}
  {( 1 - \sigma \tau )^2}
p_2(x)$$
$$+
\frac{
    {\left( 1 - \sigma  \tau  \right) }^2 +
  (\eta+\theta\sigma)(\theta+\eta\tau)}{(1-\sigma\tau)^3}{(t+\tau)}p_1(x)
    \,,$$
$$x p_n(x)=(1-\sigma\tau)^{n-1}(1+\sigma t)p_{n+1}(x)+
\frac{(\eta  + 2\theta \sigma  + \eta \sigma \tau )t+\theta  + 2\eta \tau  + \theta \sigma \tau}
  {{\left( 1 - \sigma \tau  \right) }^2}
p_{n}(x)$$
$$
+\frac{
    {\left( 1 - \sigma  \tau  \right) }^2 +
    (\eta+\theta\sigma)(\theta+\eta\tau)}
    {(1-\sigma\tau)^{n+2}}(t+\tau)p_{n-1}(x)$$
for $n\geq 2$.

After renormalizing the $n$-th polynomial by $(1-\sigma\tau)^{(n-2)(n-1)/2}$, 
from  Theorem \ref{Thm C-com+} we get the following.
\begin{proposition}[free quadratic harnesses]\label{EX free biM}
Suppose $(X_t)$ is a quadratic harness with parameters such that
$( 1 - \sigma  \tau )^2 +(\eta+\theta\sigma)(\theta+\eta\tau)> 0$, $0\leq \sigma\tau<1$, and
$q=-\sigma\tau$. If   for $t>0$ random variable
 $X_t$ has all moments and infinite support, then
it has orthogonal martingale polynomials given by the three step recurrences \rf{generic 0},  \rf{generic 1}, and
\begin{eqnarray*}
x p_2(x;t)&=&(1+\sigma t)p_3(x;t)+
\frac{(\eta  + 2\,\theta \,\sigma  + \eta \,\sigma \,\tau )t+\theta  + 2\,\eta \,\tau  + \theta \,\sigma \,\tau}
  {{\left( 1 - \sigma \,\tau  \right) }^2} p_2(x;t)
  \\&&+
 \frac{
   \left( 1 - \sigma  \tau  \right)^2 +
    +(\eta+\theta\sigma)(\theta+\eta\tau)}{(1-\sigma\tau)^3}(t+\tau)p_1(x;t), \label{Free2}\\ \nonumber
x p_n(x;t)&=&(1+\sigma t)p_{n+1}(x;t)+
\frac{(\eta  + 2\,\theta \,\sigma  + \eta \,\sigma \,\tau )t+\theta  + 2\,\eta \,\tau  + \theta \sigma \tau}
  {{( 1 - \sigma \tau) }^2}p_n(x;t)
  \\&&+
\frac{
    {\left( 1 - \sigma  \tau  \right) }^2 +(\eta+\theta\sigma)(\theta+\eta\tau)}{(1-\sigma\tau)^{4}}(t+\tau)p_{n-1}(x;t), \; n\geq 2. \label{Free n}
\end{eqnarray*}
\end{proposition}

\begin{remark}
The recurrence in Proposition \ref{EX free biM} is a finite perturbation of the constant coefficient recurrence
which was analyzed by many authors, see
\cite{Sansigre-Valent-95} and the references therein, see also
 \cite{Ronveaux-Van-Assche-95}.
\end{remark}

\begin{remark}
In \cite{Bryc-Wesolowski-04} we show that the free bi-Poisson process is associated with the generalized free convolution
studied in  \cite{Bozejko-Wysoczanski01}.
It is interesting to ask if an analogous situation occurs for the general free harnesses of Proposition
\ref{EX free biM};  for recent extensions of generalized free convolutions to perturbations of higher-order terms,
 see \cite{Krystek-Wojakowski04}.
\end{remark}
%
%
\subsection{Classical quadratic harnesses}\label{Sect classical}
The classical quadratic harnesses have parameter $q=1-2\sqrt{\sigma\tau}$.
%
 The  adjective "classical"  is motivated by the fact that when $\sigma\tau=0$, quadratic harnesses
 with $q=1$ are related to classical probability,  see \cite[Section 4.2]{Bryc-Wesolowski-03}.

\begin{proposition}[Classical quadratic harnesses]\label{EX classical biM}
Suppose $(X_t)$ is a quadratic harness with parameters such that $( 1 - \sigma  \tau )^2 +(\eta+\theta\sigma)(\theta+\eta\tau)> 0$, $0\leq \sigma\tau<1$, and
$q=1-2\sqrt{\sigma\tau}$. If   for $t>0$ random variable
 $X_t$ has all moments and infinite support, then
 it has orthogonal martingale polynomials given by the three step recurrences \rf{generic 0},  \rf{generic 1}, and
\begin{eqnarray}
x p_n(x;t)&=&\left(1+(n-1)\rho+n\sigma t \right)p_{n+1}(x;t)+
\left(\gamma_n t +\delta_n\right)p_n(x;t)
  \\&&+
\omega_n^\circ \left((1+(n-2)\rho) t +\tau (n-1)\right) p_{n-1}(x;t), \; n\geq 2. \label{Classic n}
\end{eqnarray}
where $\rho=\sqrt{\sigma\tau}$,
\begin{multline}\label{gamma c sol}
\gamma_n=\frac{n(1+(n-2)\rho)}{(1-\rho)^2(1+(2n-1)\rho)(1+(2n-3)\rho)}\times\Big(
\eta+(2n-1)\theta\sigma \\
+(2n-3)\eta\rho+2(n-1)^2\eta\rho^2+(2(n-1)^2-1)\theta\sigma\rho
\Big), \; n\geq 1
\end{multline}
\begin{multline}\label{delta c sol}
\delta_n=\frac{n(1+(n-2)\rho)}{(1-\rho)^2(1+(2n-1)\rho)(1+(2n-3)\rho)}\times\Big(
\theta+(2n-1)\eta\tau \\
+(2n-3)\theta\rho+2(n-1)^2\theta\rho^2+(2(n-1)^2-1)\eta\tau\rho
\Big), \; n\geq 1
\end{multline}

\begin{multline}\label{omega c sol}
\omega_n^\circ=
\frac{n(1+(n-3)\rho)}{(1-\rho)^2(1+(2n-2)\rho)(1+(2n-4)\rho)} \\
+\frac{n(n-1)(1+(n-2)\rho)(1+(n-3)\rho)}{(1-\rho)^4(1+(2n-2)\rho)(1+(2n-3)\rho)^2(1+(2n-4)\rho)} \\
\times\big(\left(1+(n-2)\rho\right)\theta+(n-1)\eta\tau\big)\big(\left(1+(n-2)\rho\right)\eta+(n-1)\theta\sigma\big)
,\; n\geq 2.
\end{multline}

\end{proposition}
\begin{proof} 

The assumptions of Theorem \ref{Thm C-com+} are satisfied, so the orthogonal martingale polynomials exist,
and we only need to solve the recurrences
\rf{non-linear rec}, \rf{gamma-delta rec}, \rf{w rec}, and renormalize the polynomials to simplify the final three step
recurrence \rf{Classic n}.

From  \rf{non-linear rec}   we get
$$\alpha_n=(1-\rho)^{n-2}(n-1),\; n\geq 1,$$
$$ \beta_n=(1-\rho)^{n-2}(1+(n-2)\rho), \; n\geq 2.$$
Thus $\lambda_{1,k}=\beta_{k+1}$, $k\geq 0$, and a calculation gives
\begin{equation}\label{la c}
\lambda_{n,k}=(1-\rho)^{2n+k-3}(1+(2n+k-3)\rho),\; n\geq 2
\end{equation}
$$\alpha_{n+2}\beta_n-\beta_{n+2}\alpha_n= 2 (1-\rho)^{2n-1},\; n\geq 1.
$$
Since $q+\sigma\tau=(1-\rho)^2$, equations \rf{gamma-delta rec}
simplify to
\begin{multline} \label{gamma proof c}
\gamma_{n+1}=\frac{1+(2n-1)\rho}{1+(2n+1)\rho}\gamma_n+\frac{2\sigma}{1+(2n+1)\rho}\delta_n
\\
+ \frac{(1+(2n-1)\rho+2n^2\rho^2)\eta+(2n+1+(2n^2-1)\rho)\theta\sigma}{(1-\rho)^2(1+(2n+1)\rho)},
\end{multline}
\begin{multline}  \label{delta proof c}
\delta_{n+1}=\frac{1+(2n-1)\rho}{1+(2n+1)\rho}\delta_n+\frac{2\tau}{1+(2n+1)\rho}\gamma_n
\\
+ \frac{(1+(2n-1)\rho+2n^2\rho^2)\theta+(2n+1+(2n^2-1)\rho)\eta\tau}{(1-\rho)^2(1+(2n+1)\rho)}.
\end{multline}
Formula \rf{delta1 gamma1} shows that \rf{gamma c sol} and  \rf{delta c sol} hold true for $n=1$.
Assuming  that \rf{gamma c sol} and  \rf{delta c sol} hold true for some $n\geq 1$, from equations
\rf{gamma proof c} and \rf{delta proof c} a computer-assisted calculation shows
 that the formulas hold true for $n+1$ as well.

We now show that the solution of  \rf{w rec} is
\begin{multline}\label{omega c sol raw}
\omega_n=
\frac{n(1+(n-3)\rho)}{(1-\rho)^{2n-2}(1+(2n-2)\rho)(1+(2n-4)\rho)} \\
+\frac{n(n-1)(1+(n-2)\rho)(1+(n-3)\rho)}{(1-\rho)^{2n}(1+(2n-2)\rho)(1+(2n-3)\rho)^2(1+(2n-4)\rho)} \\
\times\big(\left(1+(n-2)\rho\right)\theta+(n-1)\eta\tau\big)\big(\left(1+(n-2)\rho\right)\eta+(n-1)\theta\sigma\big)
,\; n\geq 1.
\end{multline}
Indeed, a calculation shows that the formula holds true for $n=1,2$. Suppose \rf{omega c sol raw}
 holds for some $n\geq 2$.
For $k=1$, formula \rf{la c} holds also for $n=1$; thus for $n\geq 2$  recurrence \rf{w rec} simplifies to
\begin{multline*}
\omega_{n+1}=\frac{1+(2n-4)\rho}{(1-\rho)^2(1+2n\rho)}\omega_n
+
\frac{1+\tau\gamma_n^2+\sigma\delta_n^2+\theta\gamma_n+\eta\delta_n-2\rho\delta_n\gamma_n
}
{(1-\rho)^{2n}(1+2n\rho)}.
\end{multline*}
Using \rf{gamma proof c}, \rf{delta proof c} and \rf{omega c sol raw}, a computer assisted calculation verifies that
\rf{omega c sol raw} holds true for $n+1$.

Renormalizing the $n$-th polynomial in \rf{three step abstract}
 by the factor $(1-\rho)^{(n-2)(n-1)/2}$, we get \rf{Classic n} with $\omega_n$ replaced by $\omega_n^\circ$.
\end{proof}

\subsection{Orthogonal martingale polynomials when $\sigma\tau=0$}\label{OM sigma-tau=0}
%

We use the standard $q$-notation
\begin{eqnarray*}
{[n]_{q}} &=&1+q+\dots +q^{n-1}, \\
{[n]_{q}!} &=&[1]_{q}[2]_{q}\dots [n]_{q}\,,
\end{eqnarray*}%
with the usual conventions $[0]_{q}=0,[0]_{q}!=1$.
In this notation, Remark \ref{R: Explicit solution} gives $\alpha_n=[n-1]_q$, and $\beta_n=1$.

Passing to the time-inverse $(tX_{1/t})$ if necessary, without loss of generality we may assume that $\sigma=0$.
In this case, recurrences in Theorem \ref{Thm C-com+} have explicit solutions  and orthogonal martingale polynomials are monic.

\begin{theorem}\label{Thm sigma=0} Suppose $(X_t)$ is a quadratic harness with covariance \rf{cov} and  parameters
such that
 $\sigma=0$, $-1<q\leq 1$, and $1+[n]_q\eta\theta+[n]_q^2\tau\eta^2>0$ for all $n$.
If for each $t>0$ random variable $X_t$
has infinite support,   then the monic orthogonal martingale
polynomials $p_n(x;t)$ are given by the recurrence
\begin{multline}
\label{recur sigma=0}
x p_n(x;t)=p_{n+1}(x;t) + \left(\eta t +\theta + ([n]_q+[n-1]_q)\eta\tau\right)[n]_qp_n(x;t) \\
+
(t +\tau [n-1]_q) \left(1+[n-1]_q\eta\theta+[n-1]_q^2\tau\eta^2\right)[n]_q p_{n-1}(x;t),\, n\geq 1,
\end{multline}
with the initial condition $ p_0=1, p_1(x)=x$.
\end{theorem}

\begin{remark}\label{Remark sharper bounds}
It is easy to give simple sufficient conditions on the parameters  which imply that
 $\omega_n:=1+[n-1]_q\eta\theta+[n-1]_q^2\tau\eta^2>0$ for all $n\geq 1$.  Suppose that
$-1<q\leq 1$.
\begin{enumerate}
\item If $\eta\theta\geq 0$ then, trivially, $\omega_n>0$.
\item If $\tau=0$,  then $\omega_n>0$ provided $1+\eta\theta>\max\{q,0\}$. Indeed, for $n\geq 2$ we have
\begin{eqnarray}
\label{q<0}
1+q \leq [n]_q \leq 1 &&\mbox{ if  $-1<q<0$}, \\
\label{q>0}
1 \leq [n]_q < 1/(1-q) &&\mbox{ if  $0\leq q< 1$}.
\end{eqnarray}
For $\eta\theta<0$ we get $\omega_n>0$ from the right hand sides of these equations.
\item if $-1<q<1$ and $\theta^2<4\tau$ then $\omega_n>0$. Indeed, then the quadratic function
$f(x)=1+x\eta\theta+x^2\tau\eta^2$ is non-negative.
\item If $\tau>0, \eta\theta<0$, $1+\eta\theta+\tau\eta^2>0$, and $\theta^2\geq4\tau$, then the sufficient condition is
\begin{equation}\label{+++}
1+\min\{q,0\}>\frac{|\theta|+\sqrt{\theta^2-4\tau}}{2 \tau |\eta|}.
\end{equation}
Indeed, inequality $1+\eta\theta+\tau\eta^2>0$ implies $\omega_2>0$.
Considering separately the cases $\theta<0$ and $\theta>0$,
the larger root of $1+x\eta\theta+x^2\tau\eta^2=0$ is
$\frac{|\theta|+\sqrt{\theta^2-4\tau}}{2 \tau |\eta|}$. If \rf{+++} holds then the left hand sides of inequalities
\rf{q<0}, \rf{q>0} imply $\omega_n>0$ for $n\geq 3$.
\end{enumerate}

\end{remark}
It might be interesting to point out the explicit recurrence for the case $\tau=0$. 
\begin{corollary}\label{Cor tau=0} Suppose $(X_t)$ is a quadratic harness with covariance \rf{cov} and  parameters
such that
 $\tau=0$, $-1<q<1$, and $1+[n]_q\eta\theta+[n]_q^2\sigma\theta^2>0$ for all $n$.
 Suppose that for $t>0$ random variable $X_t$ has infinite 
 support. 
Then as the orthogonal martingale
polynomials for $(X_t)$ we can take polynomials $p_n(x;t)$ given by the recurrence
$$
x p_n(x;t)=(1+t\sigma[n]_q)p_{n+1}(x;t)$$
$$+(\theta+t\eta+([n]_q+[n-1]_q)\theta\sigma t)p_n(x;t)+ t(1+[n-1]_q\eta\theta+[n-1]_q^2\sigma\theta^2)[n]_q p_{n-1}(x;t).
$$
\end{corollary}

\begin{proof}[Proof of Theorem \ref{Thm sigma=0}]
We use the notation from \rf{coeff}. The integrability assumption
of  Theorem \ref{Thm C-com+} is fulfilled by Theorem \ref{T4.1}
 and the system of recurrences in Theorem \ref{Thm C-com+} simplifies.
Since  $\alpha_n=[n-1]_q$, $\beta_n=1$, the first equation of recurrence \rf{gamma-delta rec} becomes
$$
\gamma_{n+1}=\eta+\gamma_n q
$$
which gives
\begin{equation}\label{sol gamma}
\gamma_n=\eta[n]_q.
\end{equation}
The other recurrences are solved inductively.
Suppose that
\begin{eqnarray}\label{inductive ass del}
\delta_{n-1}&=&(\theta+\eta\tau([n-1]_q+[n-2]_q))[n-1]_q\,, \\
 \label{inductive ass eps}
\omega_n&=&(1+[n-1]_q\eta\theta+[n-1]_q^2\tau\eta^2)[n]_q\,,
\end{eqnarray}
hold. The initial conditions say that both formulas are satisfied for $n=1$, and the second one
 holds true also for $n=2$.

We use \rf{sol gamma}  and \rf{inductive ass del} to compute $\delta_{n}$ from  
the second equation in \rf{gamma-delta rec}
as follows
\[ \delta_n
 =\theta+q\delta_{n-1}+q([n]_q-[n-2]_q)\tau\gamma_{n-1}+\eta\tau[n-1]_q+\eta\tau[n]_q
 =\theta(1+q[n-1]_q)
 \]
 \[+\eta\tau(q([n-1]_q+[n-2]_q)[n-1]_q+q([n]_q-[n-2]_q)[n-1]_q+[n-1]_q+[n]_q)
\]
 \[ =\theta[n]_q+\eta\tau(q([n-1]_q+[n]_q)[n-1]_q+[n-1]_q+[n]_q)
 \]
 \[=\theta[n]_q+\eta\tau([n-1]_q+[n]_q)(1+q[n-1]_q)
=\theta[n]_q+\eta\tau([n-1]_q+[n]_q)[n]_q
\]

Finally, \rf{w rec}
determines  $\omega_{n+1}$ as
$$
\omega_{n+1}
=q\omega_n+\gamma_n^2\tau
+q\gamma_n\delta_n+\gamma_n\theta+\delta_n\eta+1 -\gamma_n\delta_n.
$$
 Using the inductive assumption  \rf{inductive ass eps}, and already established
 formulas \rf{sol gamma} and  \rf{inductive ass del},
  we get
$$
\omega_{n+1}
=(1+q[n]_q)
+(1+q [n-1]_q)[n]_q\eta\theta+q[n-1]_q^2[n]_q\tau\eta^2
+\tau\eta^2[n]_q^2
$$
$$+\eta(1+q[n]_q)\delta_n
 -\eta[n]_q\delta_n
=[n+1]_q
+[n]_q^2\eta\theta+q[n-1]_q^2[n]_q\tau\eta^2
+\tau\eta^2[n]_q^2
+\eta q^{n}\delta_n
 $$
$$
=[n+1]_q+[n+1]_q[n]_q\eta\theta+\tau\eta^2(
q[n-1]_q[n]_q([n-1]_q+q^{n-1})+[n]_q^2+q^{n}[n]_q^2
).
$$
The coefficient at $\tau\eta^2$ simplifies to
$$
q[n-1]_q[n]_q^2+[n]_q^2+q^{n}[n]_q^2
=[n+1]_q[n]_q^2.
$$
Thus
 $\omega_{n+1}=(1+[n]_q\eta\theta+[n]_q^2\tau\eta^2)[n+1]_q$.

Formula \rf{recur sigma=0} comes from \rf{three step abstract} after substituting  (\ref{sol gamma}-\ref{inductive ass eps}) into \rf{coeff}.
\end{proof}

\subsection{Operator solutions} In this section we
re-derive the recurrences for some special
orthogonal martingale
 polynomials from Section \ref{OM sigma-tau=0}  by an operator approach which is related
 to Lie algebra techniques.
 This method has a 
 more  ad hoc character so we concentrate on
  two relatively simple cases only.

In the operator approach, we go back directly to Theorem \ref{T5.1}. We
re-interpret the matrices $\xxx,\yyy$ as the linear
operators acting on the formal power series in an auxiliary variable $z$.
This  identifies the martingale polynomial
 $p_n(x;t)$ with  $z^n$, where $z$ is an auxiliary variable; similar technique appeared in
 umbral calculus \cite[Ch. 1]{Rota}, in
 orthogonal polynomials \cite{IS03}, and in
 Segal-Bargmann representation \cite{Perelomov}.

We seek
 the solutions of \rf{C-com} in terms  of  the
$q$-differentiation operator
$$\D(g)(z)=\begin{cases}\frac{g(z)-g(qz)}{(1-q)z} &q\ne 1,\\
g'(z) & q=1,
\end{cases}$$
and the multiplication operator
$$
\Z (g)(z)=zg(z),
$$
treated as the linear operators on formal series $g(z)$ in variable $z$.
Table \ref{Table 1} lists the $q$-commutators of the combinations of these two operators
that we need here and in Section \ref{Coherent states}.

The requirement on the $0$-th column of $\mC_t$ reduces to the requirement that  on the unit constant function
${\xxx} 1=0$, ${\yyy} 1 =z$. When $\sigma=0$, we have $a_n(t)=1$ so
we are looking for the operators that satisfy
\begin{equation} \label{BLD}
(t \xxx+\yyy)z^n=z^{n+1}+\mbox{lower order terms}.
\end{equation}
\begin{table}[htb]
 \footnotesize
 \begin{tabular}{|r||c|c|c|c|} \hline
 &&&& \\
$_A\backslash^ B$& $\Z$ & $\Z\D$& $\Z\D^2$ &$\Z^2\D$ \\ \hline\hline
 &&&& \\
 $\D$ & $\I$ &  $\D$ & $\D^2$ & $(1+q)\Z\D-q(1-q)\Z^2\D^2$\\ &&&& \\
 $\Z\D$& $\Z$ &   $(1-q)\Z\D+q(1-q)\Z^2\D^2$  && $\Z^2\D$\\ &&&& \\
$\Z\D^2$ &  $(1+q)\Z\D-q(1-q)\Z^2\D^2$& $\Z\D^2$ & &\\ &&&& \\
 $\Z^2\D$& $\Z^2$ &&&
 \\
 \hline
\end{tabular}
\medskip

\caption{$q$-commutators $[A,B]_q:=AB-qBA$. \label{Table 1}}
\end{table}

\begin{example}[$q$-Meixner processes] \label{Example q-M}
The $q$-Meixner process is a quadratic harness  with parameters
$\sigma=\eta=0$, see \cite{Bryc-Wesolowski-03}.
Inspecting Table \ref{Table 1}  we verify that
$$
\xxx=\D
$$
and
$$\yyy=\Z(1+\theta \D+\tau \D^2)
$$
solve \rf{C-com} when $\sigma=\eta=0$. Indeed,
$$[\xxx,\yyy]_q=[\D, \Z]_q+\theta[\D,\Z\D]_q+\tau[\D,\Z\D^2]_q=
\I+\theta \D+\tau \D^2.
$$
A calculation shows now that
$$(t\xxx+\yyy)z^n=z^{n+1}+\theta [n]_qz^n+(t+\tau[n-1]_q)[n]_q z^{n-1},$$
which by the identification of $z^n$ with $p_n(x;t)$  implies that the corresponding martingale polynomials
 satisfy the three step recurrence
\begin{equation}\label{rec-G}
xp_n(x;t)=p_{n+1}(x;t)+\theta [n]_qp_n(x;t)+(t+\tau[n-1]_q)[n]_q p_{n-1}(x;t), \; n\geq 1.
\end{equation}
Of course, this is a special case of \rf{recur sigma=0} corresponding to $\eta=0$.
Feinsilver \cite[Section 3.4]{Feinsilver-90} gives a
 reparametrization of recurrence \rf{rec-G}, and considers $q$-commutator relations $[\xxx,\yyy]_q=h$
 as well as their realizations via operators
$\D,\Z$ without considering \rf{C-com}.
Anshelevich \cite[Remark 6]{Anshelevich01} discusses a reparametrization of  recurrence \rf{rec-G} in relation to free Sheffer systems.
In \cite{Bryc-Wesolowski-03} recurrence \rf{rec-G} is  used as the first step in
 the construction
of the $q$-Meixner Markov processes.
\end{example}
We now use the same method to derive the recurrence for the martingale polynomials of the bi-Poisson process.
\begin{example}[bi-Poisson process]\label{T5.2}
The bi-Poisson process is a quadratic harness  with parameters
$\sigma=\tau=0$. 
Inspecting  Table \ref{Table 1} we verify that
$$
\xxx=\D+\eta\Z(\D+\theta\D^2)
$$
and
$$\yyy=\Z(1+\theta \D)
$$
solve \rf{C-com} when $\sigma=\tau=0$.
Indeed,
$$[\xxx,\yyy]_q=[\D, \Z]_q+\theta[\D,\Z\D]_q+\eta[\Z\D,\Z]_q +\eta\theta [\Z\D,\Z\D]_q
$$
$$+
\eta\theta[\Z\D^2,\Z]_q+\eta\theta^2[\Z\D^2,\Z\D]_q
=\I+ \theta \D+\eta \Z$$
$$+
\eta\theta(1-q)\Z\D+\eta\theta q(1-q)\Z^2\D^2+
\eta\theta(1+q)\Z\D - \eta\theta q(1-q)\Z^2\D^2 + \eta\theta^2 \Z\D^2
$$
$$
=\I+ \theta \D+\eta \Z+\eta\theta  \Z\D+ \eta\theta \Z\D+\eta\theta^2\Z\D^2
$$
$$
=\I+ (\theta \D +\eta\theta \Z\D+\eta\theta^2\Z\D^2)+(\eta \Z+\eta\theta  \Z\D)
$$
$$
=\I+\theta\xxx+\eta\yyy.
$$
Operator $t\xxx+\yyy$ satisfies the constraint $(t\xxx+\yyy)1=z$. A calculation shows that
$$(t\xxx+\yyy)z^n= z^{n+1}+(\theta+t\eta)[n]_qz^{n}+t(1+\eta\theta[n-1]_q)[n]_qz^{n-1},$$
so the constraint \rf{BLD} holds.
By the identification of $z^n$ with $p_n(x;t)$, the corresponding martingale polynomials
 satisfy the three step recurrence
\begin{equation}\label{Q-Poisson}
x p_n(x;t)=p_{n+1}(x;t)+(\theta+t\eta)[n]_qp_n(x;t)+t(1+\eta\theta[n-1]_q)[n]_qp_{n-1}(x;t),
\end{equation}
$n\geq 0$, with $p_{-1}=0,p_0=1$.  Of course, this is a special case of \rf{recur sigma=0} corresponding to $\tau=0$.
But this recurrence was hard to guess without \rf{C-com}, so in \cite{Bryc-Wesolowski-04} it appears
 for  $q=0$ only.

\end{example}

\subsection{Dual $q$-commutation equation. Coherent states}\label{Coherent states}
Coherent states and the Segal-Bargmann representation are  analytical techniques developed in mathematical physics
\cite{Perelomov}.
The full Segal-Bargmann isomorphism is known to fail
even in the relatively simple case of $q$-Brownian motion with $q<0$, compare
\cite{Maassen95a}.
But  algebraic duality is available and useful in a more general setting.

Let
\begin{equation}\label{NC Q*}
Q_{t,s,u}^*(\xxx,\yyy)=\Atsu \xxx^2+\Btsu \yyy\xxx +\Ctsu \yyy^2+\Dtsu \xxx + \Etsu \yyy +\Ftsu
\end{equation}
be the  quadratic
form  in the non-commuting variables $\xxx,\yyy$; for background, see e.g. \cite[page 7]{Kassel95}.
Note that this is a dual of the quadratic form \rf{NC Q} that appears in the proof of Theorem \ref{T5.1}.
The following  relates  \rf{NC Q*} to the dual version of the $q$-commutation equation \rf{C-com}.

\begin{proposition}\label{P5.2}
Let  $\X_t$ be a linear function $\X_t=\xxx+t \yyy$ of non-commutative variables $\xxx, \yyy$, and suppose that
the coefficients of the
quadratic form $Q_{t,s,u}^*$ are given by \rf{F0} and (\ref{A0}-\ref{E0}).
Then the following statements are equivalent.
\begin{enumerate}
\item Operator identity
\begin{equation}\label{Eqn: operator QH}
\X_t^2=Q_{t,s,u}^*(\X_s,\X_u)
\end{equation}
holds for all $s<t<u$;
\item The non-commutative variables $\xxx, \yyy$  satisfy the
equation
\begin{equation}\label{EQ: q-com}
[\xxx,\yyy]_q=\sigma \xxx^2+\tau \yyy^2+\eta \xxx+\theta \yyy + 1.
\end{equation}
\end{enumerate}
\end{proposition}
\begin{proof} 
Since $\X_u\X_s=\xxx^2+su \yyy^2+ s(\xxx\yyy+\yyy\xxx)+(u-s)\yyy\xxx$, we have
\begin{multline*}\X_t^2-Q_{t,s,u}^*(\X_s,\X_u)\\
=
(1-\Atsu-\Btsu-\Ctsu)\xxx^2+ (t^2-s^2\Atsu-su\Btsu-u^2\Ctsu) \yyy^2\\
+ (t-s\Atsu-s\Btsu-u\Ctsu)(\xxx\yyy+\yyy\xxx)-(u-s)\Btsu\yyy\xxx
\\
-(\Dtsu+\Etsu)\xxx-(s\Dtsu+u\Etsu)\yyy-\Ftsu\I.
\end{multline*}
Relations \rf{0}, \rf{sigma-tau}, \rf{q-eta-theta} are equivalent to
\[\X_t^2-Q_{t,s,u}^*(\X_s,\X_u)=
-\sigma\Ftsu\xxx^2 - \tau\Ftsu \yyy^2 \]
\[+ \Ftsu(\xxx\yyy+\yyy\xxx)-(1+q)\Ftsu\yyy\xxx
-\eta\Ftsu\xxx-\theta\Ftsu\yyy-\Ftsu\I.
\]
\[=\Ftsu([\xxx,\yyy]_q-\sigma\xxx^2-\tau\yyy^2-\eta\xxx-\theta\yyy-\I)=0.
\]
\end{proof}

\subsubsection{Coherent states of  $q$-Meixner process}\label{Appl: Meixner}
Continuing Example \ref{Example q-M}, the analogue of the coherent state in physics is the generating function
\begin{equation}\label{coherent}
\varphi_t(z,x)=\sum_{n=0}^\infty \frac{z^n}{[n]_q!}p_n(x;t),
\end{equation}
where $p_n(x;t)$ are the  orthogonal martingale polynomials given by recurrence
\rf{rec-G}.

\begin{proposition} \label{Prop QM}The operator
\begin{equation}\label{Eqn: XX=}
\X_t=(1+\theta \Z+\tau \Z\;^2)\D+t \Z
\end{equation}
satisfies \rf{Eqn: operator QH} with the
quadratic form $Q_{t,s,u}^*$ given by \rf{F0}, (\ref{A0}-\ref{E0}) where $\sigma=\eta=0$. Moreover,
\begin{equation}\label{SB}
x\varphi_t(z,x) = (\X_t\varphi_t)(z,x).
\end{equation}
\end{proposition}
\begin{proof}
Inspecting  Table \ref {Table 1}  it is easy to verify that
 $\xxx=(1+\theta \Z+\tau \Z^2)\D$ and $\yyy=\Z$
 satisfy \rf{EQ: q-com} with $\sigma=\eta=0$.
 Indeed,
 \[
 [\xxx,\yyy]_q=[\D,\Z]_q+\theta[\Z\D,\Z]_q+\tau[\Z^2\D,\Z]_q
 \]
 \[
 =\I+\theta \Z+\tau \Z^2=\I+\theta\yyy+\tau\yyy^2.
 \]
 Therefore,  Proposition \ref{P5.2} implies \rf{Eqn: operator QH}.
The algebraic identity \rf{SB} follows from recurrence \rf{rec-G} by the following calculation.
\[
x\varphi_t(z,x)=\sum_{n=0}^\infty \frac{z^n}{[n]_q!}xp_n(x;t)
\]
\[=
\sum_{n=0}^\infty \frac{z^n}{[n]_q!}\left(p_{n+1}(x;t)+\theta [n] p_n(x;t)+\tau[n][n-1]p_{n-1}(x;t)+t [n]p_{n-1}(x;t)\right)=
\]
\[=
\D \sum_{n=0}^\infty \frac{z^{n+1}}{[n+1]_q!}p_{n+1}(x;t)+\theta z \D \sum_{n=0}^\infty \frac{z^n}{[n]_q!}p_{n}(x;t)
\]
\[+
\tau z^2 \D\sum_{n=1}^\infty \frac{z^{n-1}}{[n-1]_q!}p_{n-1}(x;t)+tz \sum_{n=1}^\infty \frac{z^{n-1}}{[n-1]_q!}p_{n-1}(x;t)
\]
\[
=(\D+\theta\Z\D+\tau \Z^2\D+t\Z)\varphi_t(z,x).
\]

\end{proof}

Formula \rf{SB} is implicit in the usual derivation of the product formula for $\varphi_t(z,x)$, compare \cite{Al-Salam76}.
When the parameters $\tau, \theta$ vanish, it  appears in
\cite[Section 3]{Maassen95a} in the context of analyzing ground states for the $q$-deformed Gauss distribution.
Ref. \cite{Asai-Kubo-Kuo-03} gives a more general analytical scheme which coincides with \rf{Eqn: XX=} and \rf{SB} when $\tau=0,q=1$;
however, it advocates the normalization by the $L_2$-norm of the polynomials,
which does not fit all the cases
we are interested in.

In \cite{Bryc-Wesolowski-03} we defined  the $q$-Meixner process as a Markov process with the initial state $X_0=0$ and
 with the transition probabilities $P_{s,t}(x,dy)$
determined as the unique probability measure orthogonalizing  the
polynomials $Q_n$ in variable $y$ which are given by the
three step recurrence \begin{eqnarray}\label{Q-rec-G}  && yQ_n(y|x) \\ \nonumber
&&=
Q_{n+1}(y|x,t,s) +(\theta [n]_q+x
q^n)Q_n(y|x,t,s)\\ \nonumber
&&+(t-sq^{n-1}+\tau[n-1]_q)[n]_q Q_{n-1}(y|x,t,s).
\end{eqnarray}
In that paper, we showed that this Markov process is well defined, that it has orthogonal martingale polynomials
 $p_n(x;t)$ given by recurrence \rf{rec-G}, and we used this to prove that $(X_t)$ is a quadratic harness
 with parameters $\sigma=\eta=0$.

Proposition \ref{Prop QM} simplifies
 the verification of the quadratic harness condition in \cite[Proposition 3.4]{Bryc-Wesolowski-03}, condensing more than
 three pages of proof into one page.
\begin{proposition}[\cite{Bryc-Wesolowski-03}]\label{Bryc-Wesolowski-03}
If $|q|\leq 1$ and the polynomials $p_n(x;t)$ defined  by \rf{rec-G} are  orthogonal martingale polynomials for a
Markov process $(X_t)$, then
$(X_t)$ is a quadratic harness with parameters $\sigma=\eta=0$.
\end{proposition}

\begin{proof} For simplicity, we consider only $|q|<1$, as in this case \rf{rec-G} implies that
$|X_t|\leq C t$ has bounded support.
Consider the generating function \rf{coherent}.
Since $(X_t)$ is Markov,
condition \rf{LH} is equivalent to
\begin{eqnarray}\lbl{Eqn: LH-chi}
&&\wo\left(\varphi_s(z_1,X_s)X_t\varphi_u(z, X_u)\right)\\ \nonumber
&&=\aaa \wo\left(\varphi_s(z_1,X_s)X_s\varphi_s(z, X_s)\right)+
\bbb \wo\left(\varphi_s(z_1,X_s)X_u\varphi_u(z, X_u)\right)
\end{eqnarray}
holding for all $z_1,z$ in a neighborhood of $0$ (or just as an identity in formal power series in variables $z_1,z$).

We now use the martingale polynomial property which implies that
\begin{equation}\label{Eqn: Appell}
\wwo{\varphi_t(z,X_t)}{X_s}=\varphi_s(z,X_s),
\end{equation}
and we use \rf{SB}
to represent the process through the operator \rf{Eqn: XX=} as
$$
X_t\varphi_t(z,X_t)=\X_t(\varphi_t(z,X_t)).
$$
This gives
$$
\wo\left(\varphi_s(z_1,X_s)X_t\varphi_u(z, X_u)\right)=\wo\left(\varphi_s(z_1,X_s)\X_t\varphi_t(z, X_t)\right)
$$
$$
=\X_t\wo\left(\varphi_s(z_1,X_s)\varphi_s(z, X_s)\right)=\X_tG_s(z_1,z),
$$
where
$$G_s(z_1,z)= \wo(\varphi_s(z_1,X_s)\varphi_s(z,X_s))=
\sum_{n=0}^\infty \frac{(z_1z)^n}{[n]_q!^2}\wo(p_n(X_s;s)^2)$$
$$=\sum_{n=0}^\infty \frac{(z_1z)^n}{[n]_q!}\prod_{k=1}^n(t+\tau[k-1]),$$
and $\X_t$ acts on  $G_s(z_1,z)$ as a series in variable $z$.
Therefore, equation \rf{Eqn: LH-chi}
is equivalent to
$$
\X_tG_s(z_1,z)=\aaa \X_sG_s(z_1,z)+\bbb \X_uG_s(z_1,z),
$$
and
 follows from
 the (trivial) operator identity
$\X_t=\aaa \X_s+\bbb \X_u.$

Similarly, condition \rf{QH} is equivalent to
\begin{eqnarray}\label{Eqn: QH-chi}
&&\wo\left(\varphi_s(z_1,X_s)X_t^2\varphi_u(z_2, X_u)\right)
=
\AAA \wo\left(\varphi_s(z_1,X_s)X_s^2\varphi_s(z_2, X_s)\right) \\  &&\nonumber
+
\BBB \wo\left(\varphi_s(z_1,X_s)X_sX_u\varphi_u(z_2, X_u)\right)
+\CCC \wo\left(\varphi_s(z_1,X_s)X_u^2\varphi_u(z_2, X_u)\right)\\ &&\nonumber+
\Dtsu \wo\left(\varphi_s(z_1,X_s)X_s\varphi_s(z_2, X_s)\right) +
\Etsu \wo\left(\varphi_s(z_1,X_s)X_u\varphi_u(z_2, X_u)\right) \\ &&\nonumber
+\Ftsu \wo\left(\varphi_s(z_1,X_s)\varphi_s(z_2, X_s)\right) . \nonumber
\end{eqnarray}
Notice that for $s\leq u$ we have
$$
\wo\left(\varphi_s(z_1,X_s)X_sX_u\varphi_u(z, X_u)\right)=\wo\left(\varphi_s(z_1,X_s)X_s\X_u\varphi_u(z, X_u)\right)
$$
$$
=\X_u\wo\left(\varphi_s(z_1,X_s)X_s\varphi_s(z, X_s)\right)=
\X_u\X_s \wo\left(\varphi_s(z_1,X_s)\varphi_s(z, X_s)\right)=\X_u\X_sG_s(z_1,z).
$$
Therefore, equation \rf{Eqn: QH-chi} follows from the operator identity \rf{Eqn: operator QH}
which in expanded form says
$$\X_t^2=\AAA \X_s^2+\BBB\X_u\X_s+\CCC\X_u^2+\DDD \X_s+\EEE\X_u+\FFF \I,$$
applied to  $G_s(z_1,z)$ treated as a series in variable $z$.
\end{proof}
\subsubsection{Coherent states of the bi-Poisson process} We now repeat the methods of Section \ref{Appl: Meixner} to
derive  new results about the bi-Poisson process from Example \ref{T5.2}.
\begin{proposition}\label{Prop 3.10}The operator
\begin{equation}\label{Eqn: XXx=}
\X_t=(1+(\theta+\eta t)\Z+t\eta\theta\Z\,^2)\D+t \Z,
\end{equation}
satisfies \rf{Eqn: operator QH} for the quadratic form $Q_{t,s,u}^*$
with the coefficients  \rf{F0}, (\ref{A0}-\ref{E0}) such that $\sigma=\tau=0$. Moreover,
if $\varphi_t(z,x)$ is the generating function
\rf{coherent} of the  orthogonal martingale polynomials  $p_n(x;t)$ given by recurrence
\rf{Q-Poisson}, then
\begin{equation}\label{SB P}
x\varphi_t(z,x) = (\X_t\varphi_t)(z,x).
\end{equation}
\end{proposition}

\begin{proof}
Inspecting  Table \ref{Table 1}, we verify that
 $$\xxx=(1+\theta\Z)\D$$ and $$\yyy=\Z+\eta( \Z+\theta\Z^2)\D$$
 solve \rf{EQ: q-com} with $\sigma=\tau=0$.
 Indeed,
 \[
 [\xxx,\yyy]_q=[\D,\Z]_q+\theta[\Z\D,\Z]_q+\eta[\D,\Z\D]_q+\eta\theta[\Z\D,\Z\D]_q+
 \eta\theta [\D,\Z^2\D]_q\]
 \[+\eta\theta^2[\Z\D,\Z^2\D]_q
 =\I+\theta \Z+\eta\D+\eta\theta ((1-q)\Z\D+q(1-q)\Z^2\D^2)
\]
 \[ +
 \eta\theta ((1+q)\Z\D-q(1-q)\Z^2\D^2)+\eta\theta^2\Z^2\D
 =\I+\theta \Z+\eta\D+\eta\theta \Z\D
\]
 \[ +
 \eta\theta \Z\D+\eta\theta^2\Z^2\D
 =\I+\eta(\I+\theta \Z)\D
 +
 \theta (\Z+\eta( \Z\theta \Z^2)\D)=\I+\eta\xxx
 +
 \theta \yyy.
 \]
 Therefore,  Proposition \ref{P5.2} implies \rf{Eqn: operator QH}.

 We now derive    \rf{SB P}   from \rf{Q-Poisson} by the following calculation.
\[
x\varphi_t(z,x)=\sum_{n=0}^\infty \frac{z^n}{[n]_q!}xp_n(x;t)
\]
\[=
\sum_{n=0}^\infty \frac{z^n}{[n]_q!}\left(
p_{n+1}(x;t)+(\theta+t\eta)[n]_qp_n(x;t)+t(1+\eta\theta[n-1]_q)[n]_qp_{n-1}(x;t)\right)=
\]
\[=
\D \sum_{n=0}^\infty \frac{z^{n+1}}{[n+1]_q!}p_{n+1}(x;t)+(\theta + t\eta) z \D \sum_{n=0}^\infty \frac{z^n}{[n]_q!}p_{n}(x;t)
\]
\[+tz \sum_{n=1}^\infty \frac{z^{n-1}}{[n-1]_q!}p_{n-1}(x;t) +
t\eta\theta z^2 \D\sum_{n=1}^\infty \frac{z^{n-1}}{[n-1]_q!}p_{n-1}(x;t)
\]
\[
 =\left((1+(\theta+\eta t)\Z+t\eta\theta\Z^2)\D+t \Z\right)
\varphi_t(z,x).
\]

\end{proof}
Next we give  the bi-Poisson version of Proposition \ref{Bryc-Wesolowski-03}.
\begin{proposition}\label{Prob biPoisson}
Fix  $-1\leq q\leq 1$.
If polynomials $p_n(x;t)$ given by \rf{Q-Poisson} are
orthogonal martingale polynomials for a
Markov process $(X_t)$, then
$(X_t)$ is a bi-Poisson process, i.e.,
\rf{LH} and
\rf{q-Var} hold with $\sigma=\tau=0$, and $1+\eta\theta\geq \max\{q,0\}$.
Moreover, such Markov process $(X_t)$ is determined uniquely.
\end{proposition}

\begin{proof} 
Property \rf{cov} follows from the explicit form of the polynomials $p_0,p_1,p_2$.
Orthogonality gives $\wo(X_t)=\wo(p_1(X_t;t)p_0(X_t;t))=0$ and
$\wo(X_t^2)=\wo((p_2(X_t;t)+(\theta+\eta t)X_t+t)p_0(X_t))=t$. Martingale polynomial property then implies
 $\wo(X_sX_t)=\wo(X_s\wwo{p_1(X_t;t)}{\cal F_{\leq s}})=\wo(X_s^2)=s$.

If $-1\leq q <1$, from \rf{Q-Poisson} it follows that$|X_t|\leq C t$  has bounded support;
if $q=1$ then it is not hard to identify the distribution
of $|X_t|$ for example, from \cite[Chapter 4]{Schoutens00} or \cite[pages 175--181]{Chihara}, verifying that   $|X_t|$
has a finite exponential moment. Thus the moment problem has unique solution and
the process $(X_t)$ is determined uniquely. Moreover,
 polynomials are dense in $L_2(X_s,X_u)$, see \cite[Theorem
3.1.18]{Dunkl-Xu}.

Consider the generating function \rf{coherent} with polynomials $p_n(x;t)$
 satisfying \rf{Q-Poisson}. Notice that property  \rf{Eqn: Appell} follows again from
 the  martingale polynomial condition.
 By
 Proposition \ref{Prop 3.10} with $\X_t$ defined by \rf{Eqn: XXx=}   we have the representation
$$
X_t\varphi_t(z,X_t)=\X_t(\varphi_t(z,X_t)).
$$

Since $(X_t)$ is Markov and  polynomials are dense in $L_2(X_s,X_u)$,
condition \rf{LH} is again equivalent to \rf{Eqn: LH-chi}, which we can
interpret
 as the identity between the formal power series in variables $z_1,z$. The latter follows again from
 the (trivial) operator identity $\X_t=\aaa \X_s+\bbb \X_u$, applied to
$$G_s(z_1,z)= \wo(\varphi_s(z_1,X_s)\varphi_s(z,X_s))=
\sum_{n=0}^\infty \frac{(z_1z)^n}{[n]_q!^2}\wo(p_n(X_s;s)^2),$$
treated as the formal power series in variable $z$.

Similarly, condition \rf{QH} is equivalent to \rf{Eqn: QH-chi}, and again
$$
\wo\left(\varphi_s(z_1,X_s)X_sX_u\varphi_u(z, X_u)\right)=\X_u\X_sG_s(z_1,z).
$$
Therefore, equation \rf{Eqn: QH-chi} follows from the operator identity \rf{Eqn: operator QH}
applied to  $G_s(z_1,z)$ treated as the formal power series in variable $z$.

Since the third coefficient in \rf{Q-Poisson} is nonnegative for all $n$, we get $1+\eta\theta\geq \max\{q,0\}$,
compare Remark \ref{Remark sharper bounds}.
\end{proof}

\subsection*{Acknowledgement} We would like to thank
M. Bo\.zejko for numerous discussions and hospitality,
L. Gallardo for preprint of  \cite{Gallardo-Yor04b},
 M. Ismail, D. Pommeret,   M. Yor and V. Zarikjan for helpful discussions.
The second and the third authors are very grateful to their hosts for
providing excellent research facilities and friendly atmosphere during
their visit to the Department of Mathematics, University of Cincinnati,
in August and September 2004.


\end{document}